\documentclass[12pt]{article}
\usepackage{setspace}

\usepackage{amsmath}
\usepackage{amssymb}
\usepackage{multirow}
\usepackage{color}
\renewcommand{\epsilon}{\varepsilon}

\newcommand{\N}{\mathbb N}
\newcommand{\R}{\mathbb R}

\usepackage{amsmath}
\usepackage{amssymb}
\usepackage[ansinew]{inputenc}
\usepackage{graphicx}
\usepackage{epsfig}

\usepackage[authoryear]{natbib}
\newcommand{\be}{\begin{eqnarray}}
\newcommand{\ee}{\end{eqnarray}}
\newcommand{\bea}{\begin{eqnarray*}}
\newcommand{\eea}{\end{eqnarray*}}

\newcommand{\ve}{\varepsilon}
\newcommand{\Real}{\mathbb R}

\newtheorem{theorem}{Theorem}[section]
\newtheorem{lemma}{Lemma}[section]

\newtheorem{corollary}{Corollary}[section]

\newtheorem{remark}{Remark}[section]
\newcommand{\HD}{\textcolor{red}}


\title{Reproducing kernel Hilbert spaces, {polynomials} and the classical moment problem\HD{}}
\begin{document}

\author{ Holger Dette   \\
 Ruhr-Universit\"at Bochum,
 Fakult\"at f\"ur Mathematik,
 44780 Bochum,  Germany
\and
 Anatoly Zhigljavsky \\
 School of Mathematics,
  Cardiff  University,
   Cardiff, CF24 4AG,   UK
}
\date{}
\maketitle
\begin{abstract}
{\small
We show  that polynomials do not belong to the reproducing kernel Hilbert space of  infinitely differentiable translation-invariant kernels
whose spectral measures  have moments corresponding to a determinate moment problem. Our proof is based on  relating this question to the problem
of the best linear estimation in continuous time one-parameter regression  models with a stationary error process defined by the kernel. In particular, we show that the existence of a sequence of estimators  with variances converging to $0$ implies that the regression function cannot be an element of the reproducing kernel Hilbert space. This question  is  then
related to the determinacy  of  the Hamburger moment problem for the spectral measure corresponding to the kernel.

%In the literature it was observed that a non-vanishing  constant function does not belong to the  reproducing kernel Hilbert space associated with the Gaussian kernel
%\citep[see Corollary 4.44 in][]{steinwart2008support}.  Our results provide a unifying view of this phenomenon and {show that the mentioned
% result can be extended for arbitrary polynomials and } a broad class of
%translation-invariant kernels.
%Moreover, we also demonstrate that similar observations can be made for general polynomials.
}
\end{abstract}

\medskip



%\begin{keyword}
%{linear regression},
%{correlated observations},
%{signed measures},
%{optimal design},
%{BLUE},
%\end{keyword}

% \end{frontmatter}

%\tableofcontents

%\newpage

\section{Introduction} \label{sec1}
\def\theequation{1.\arabic{equation}}
\setcounter{equation}{0}

\subsection{Main results}
\label{sec:main}

Let $X \! \subseteq \!\mathbb{R}^d$, $d \!\geq \!1$,
 $K\!:\! X\times X \!\to \!\mathbb{R}$ be a  positive definite kernel on $X$
  and define  $H(K)$ as  the corresponding Reproducing Kernel Hilbert Space (RKHS).
We assume that {$X$ has a non-empty interior} and $K$ is an infinitely differentiable (on the diagonal) translation-invariant kernel so that $K(x,y)\!=\!k(x\!-\!y)$, where $k: \mathbb{R}^d \to \mathbb{R}$
is a non-constant positive definite  function infinitely differentiable at the point $0$. {Clearly, the function $k(\cdot)$ is  even so that $k(u)=k(-u)$ for all  $u \in \mathbb{R}^d$.} Without loss of generality, we suppose $k(0)=1$.

\medskip

It is well-known, see e.g. Corollary 4.44 in \cite{steinwart2008support}, that in the case of the Gaussian (squared exponential)  kernel
%,where
\be \label{eq:gauss}
k(x)=\exp\{-\lambda \|x\|^2\}\;\;{\rm  with}\; \lambda>0,\ee
 the constant function does not belong to $H(K)$. {This result has been generalized to arbitrary polynomials in \cite{minh2010some}.
 The purpose of this paper is to significantly extend these results (previously only  known  for the case of the Gaussian  kernel) to  a substantially  larger class of kernels.}
\medskip

%
% \begin{itemize}
% \item[(a)]  a substantially  larger class of kernels, for which this statement remains true,  and
% \item[(b)]  arbitrary polynomials rather than simply constant functions.
% \end{itemize}

 %The results of this paper have definite consequences for the methodologies of
%function approximation,  Bayesian  global optimization  and support vector machines (SVM) and other kernel-based  machine learning, see Section~\ref{sec:int_2}.

\medskip

In the main part of the paper, we consider the case  $X \subset \mathbb{R}$. In this case, by Bochner's theorem
\citep{bochner1949}, there exists a measure $\alpha$, such that
the function $k$ can be represented in the form
\be \label{eq:spectral}
k(x)= \int_{-\infty}^{\infty} e^{i t x} \alpha(dt) ~~\mbox{ for all  } x \in X.
\ee
The measure $\alpha$ is called {\it  spectral measure}.   {As function $k(\cdot)$ is even  and $k(0)=1$, $\alpha(dt)$  is a probability measure  symmetric around the point 0.}
 {The moments of this measure (in the case of their existence) will be  denoted  by
%\be \label{eq:spectral1}
$c_n= \int_{-\infty}^{\infty} t^n \alpha(dt)$. Since we assume $k(\cdot)$  in \eqref{eq:spectral}
is infinitely differentiable and even, we have
%In the case of existence relation \eqref{eq:spectral}  and  symmetry of $k(\cdot)$ yield
\be \label{eq:spectral1}
c_n=   \int_{-\infty}^{\infty} t^n \alpha(dt)=   \left\{
                                                   \begin{array}{cl}
                                                   (-1)^{n/2} \,k^{(n)}(0), & \mbox { if $n$ is even} \\
                                                   0, & \mbox { if $n$ is odd} \\
                                                   \end{array}
                                                 \right.   \, ,
                                                    \ee
where  $k^{(n)}(0)=\frac { \partial^{n}} {\partial u^{n} } k(u)\big|_{u=0}$.
%As we have assumed that $k(\cdot)$ is infinitely differentiable at the point ~$0$ the
%moments $c_n$ exist for all $n \in \mathbb{N}$.
}

The classical {\it Hamburger moment problem}
is to give necessary and sufficient conditions such that a given real sequence $(c_{n})_{n \in \N}$ is in fact a sequence of moments of
a distribution $\alpha$   defined on the Borel sets of $\mathbb{R}=(-\infty,\infty)$. In particular, the sequence
$(c_{n})_{n \in \N}$ is  a sequence of moments  of some distribution if and only if the  {\it Hankel} matrices  $(c_{i+j})_{i,j=0, \ldots , n}$
are positive semidefinite for all $n \in \mathbb{N}$; see e.g.
  \cite{shohat1943problem,schmudgen2017moment} among many others.
The  Hamburger moment problem  is called {\it determinate}   if the sequence of moments $(c_{k})_{k \in \N}$  determines the measure $\alpha(dt)$ uniquely.

\smallskip
{
The main results of this paper are Theorems \ref{th:one_dim} and \ref{th:one_dim1} formulated below. These theorems
provide sufficient   conditions  ensuring  that the polynomials  do not belong to the RKHS
 $H(K)$.
The proofs are given  in Section~\ref{sec:proof1}. % and~\ref{sec:proof2} respectively.
}

\begin{theorem}
\label{th:one_dim}
Let  $X \subset \mathbb{R}$ and assume that the spectral measure $\alpha(dt)$ in~\eqref{eq:spectral} has infinite support
 and no mass at the point  0.
% and has infinite support.
If the  Hamburger  moment problem for this measure is determinate, then
the non-zero constant functions
%$f(x)={\rm const}\neq 0$, $\forall x \in X$,
do not belong to  the RKHS
 %Reproducing Kernel Hilbert Space
  $H(K)$.
%Moreover,  if for any integral $m $ the  Hamburger  moment problem for the measure $\alpha_m(dt) =t^{2m} \alpha(dt)/c_{2m}$  is determinate, then $H(K)$ does not contain any  non-trivial (non-zero) polynomial on $X$.
\end{theorem}

\begin{theorem}
\label{th:one_dim1}
%Let $f(\cdot)$ be a non-constant polynomial on $X$.
Let  $X \subset \mathbb{R}$, $m$ be  a positive integer and assume that the spectral measure $\alpha(dt)$ in~\eqref{eq:spectral} has infinite support. If  the  Hamburger  moment problem for the measure $\alpha_m(dt) =t^{2m} \alpha(dt)/c_{2m}$
 is determinate  {{\rm (}here, as  in \eqref{eq:spectral1}, $c_{2m}=\int t^{2m} \alpha(dt)$\rm{)}}, then the
RKHS
 %Reproducing Kernel Hilbert Space
  $  H(K)$ does not contain   polynomials on $X$ of degree  {precisely} $  m$.
  %\footnote{ {Not only polynomials but also any function whose restrictions on any subset of $X$ with a non-empty interior  is a polynomial. }}
\end{theorem}

{Theorem~\ref{th:one_dim1} can be considered as a corollary of Theorem~\ref{th:one_dim} and  therefore  the result of Theorem~\ref{th:one_dim} is more fundamental. In the very important particular case (see the sufficient condition for moment determinacy \eqref{eq:Chow} and related discussion in Section~\ref{sec:suff}), when the function $k(\cdot)$ is real analytic  and vanishes at infinity while $X$ is bounded, $H(K)$ contains only those analytic functions that vanish at infinity. This is formally proven  in \cite{karvonen2021non} and basically follows from the fact that if $f\in H(K)$ then $f(x)$ is a point-wise limit of the sums $\mu_N(x)=\sum_{i=1}^{N}w_i k(x_i-x)$ with $x_i \in X$, which necessarily vanish at infinity.
 }

Combining   Theorems~\ref{th:one_dim} and~\ref{th:one_dim1} with  their variations in the cases when the spectral measure $\alpha(dt)$ has finite support
 (see Section~\ref{sec:finite_support}) and when this measure has positive mass at 0 (see   Theorem~\ref{th:3}), we obtain the following corollary.

\begin{corollary}
\label{cor:1}
Let  $X \subset \mathbb{R}$ and  the  Hamburger  moment problem for the spectral measure $\alpha(dt) $  be determinate. Then we have the following:
\begin{itemize}
  \item[(a)] the constant functions
$f(x)={\rm const}\neq 0$, $\forall x \in X$,  belong to
 %Reproducing Kernel Hilbert Space
  $H(K)$ if and only if  $\alpha(dt)$ has a positive mass at  the point 0;
  \item[(b)]  {if  additionally $k(\cdot)$ is a real analytic function, then
 %Reproducing Kernel Hilbert Space
  $  H(K)$ does not contain non-constant  polynomials on $X$.}
\end{itemize}
\end{corollary}

Theorems \ref{th:one_dim}, \ref{th:one_dim1} and Corollary \ref{cor:1} can be easily extended to the multivariate case, see
Section~\ref{sec:multivar}.

\subsection{Implications and related results}
\label{sec:int_2}

{Gaussian process (GP) models deal with an unknown deterministic function
assuming that it is a realization of Gaussian process (field) with some mean and  covariance kernel, which are perhaps parameterized.
The popularity of the GP model comes from
its transparency, flexibility and computational tractability; it is used as a general-purpose technique to model, explore and exploit unknown functions. As a result, methods based on the GP models constitute much of the
modern statistical toolkit for function approximation, interpolation and prediction  \citep{stein1999interpolation,wendland2004scattered}, integration \citep{briol2019probabilistic},  machine learning   \citep{rasmussen2006gaussian,steinwart2008support},  space-filling \citep{pronzato2020bayesian}, signal processing \citep{cambanis1983sampling}, probabilistic numerics \citep{hennig2015probabilistic} and global optimization  \citep{zhigljavsky2021bayesian}.
The practical application areas of the GP model are vast, and we refer the  references in the cited papers  and to
the work of \cite{schulz2018tutorial}, \cite{archetti2019bayesian} among others.
}

{Consider  the common framework of   GP regression (simple kriging), where
 the  a function $f : X \to \mathbb{R} $ to be approximated ($X \subset \mathbb{R}^d$)  is considered as
  a realization of a GP, say  $\{Z_x \}_{x \in X}$,   with mean   zero, covariance
\be \mbox{$\mathbb{E}\{Z_x Z_{x'}\}=\sigma^2\, K(x,x')$ for all $x, x'\in X \subset \mathbb{R}^d$,} \label{eq:covar}
\ee
and $\sigma^2>0$ may be unknown; see
 \cite{rasmussen2006gaussian} for more details.
Let the kernel $K$ be strictly positive definite, $X_N=\{x_1,\ldots,x_N\}$ be an $N$-point design consisting of distinct points $x_i \in X$ and $F_N=[f(x_1),\ldots,f(x_N)]^\top \in \mathbb{R}^{N} $ be the vector
of exact observations of $f$ at the points of $X_N$.
%Assume that the matrix $\{K_N\}_{i,j}=K(x_i,x_j)$, $i,j=1,\ldots,n$ has full rank.
The conditional process $\{Z_x|(X_N,F_N)\}_{x \in X} $ is again Gaussian
with conditional mean
\begin{equation}
\label{det1}
\mu_N(x)=F_N^\top K_N^{-1} b_N(x)
\end{equation}
and covariance function
\be \label{det11}
C_N(x,y)=K(x,x')- b_N^\top(x) K_N^{-1} b_N(x'),
\ee
 where
$K_N=(K(x_i,x_j))_{i,j=1}^N\, $
and
$b_N(\cdot)=[K(x_1,\cdot),\ldots,K(x_N,\cdot)]^\top$. Straightforward calculation shows that the conditional mean $\mu_N(x)$ is the best linear predictor  of $f(x)$ and $\sigma^2\, C_N(x,x)$ is the corresponding  mean squared prediction error  at the point $x$.
}

{
GP regression is equivalent to kernel   interpolation, see e.g. \cite{scheuerer2013interpolation}, Section~3.3 in \cite{kanagawa2018gaussian} and Chapter 3 in \cite{paulsen2016introduction}.
More precisely,
 the conditional mean $\mu_N(\cdot)$ is the  minimal-norm interpolant to $f$  among all functions in $H(K)$, where we use the norm on $H(K)$ denoted below by $\|\cdot \|_{H(K)} $.
  This property implies in particular that if $f \in H(K)$ then $\|\mu_N \|_{H(K)} \leq \|f \|_{H(K)} $, where
   this inequality holds for any $f \in H(K)$ and any  set  of points $X_N$. If $f \notin H(K)$ then
  the conditional mean $\mu_N(\cdot)$ is still an element of $H(K)$  but its norm tends to infinity as $N$ grows and $X_N$ becomes denser.
  Therefore, there is a fundamental difference between the complexity of the  approximation problem of $f$ depending on whether  $f \in H(K)$
  or $f \notin H(K)$. Correspondingly,
 properties of all other techniques based on the use of the GP model also heavily depend on whether an unknown function of interest belongs to the corresponding
RKHS. We also refer to the work of  \cite{steinwart2006explicit} and to  Section 4.4 in \cite{steinwart2008support} for a discussion on the importance of this issue for the learning performance of
support vector machines (SVMs)  in the case of the Gaussian  kernel \eqref{eq:gauss} as well as for  the difficulty of deciding whether  a given function $f$ belongs to the RKHS    $ H(K)$ for a chosen kernel $K$. In  GP regression, knowing that the non-zero constant functions
do not belong to  the RKHS $H(K)$ is especially important as it can be used to justify omitting function centering;  see, for example, Assumption 2 in
\cite{lee2016variable}.
}

{
Assume now that the factor $\sigma^2$ in \eqref{eq:covar}  is unknown and that the maximum likelihood estimator  (MLE)
$\widehat{\sigma^2_N }$  of $\sigma^2$ is constructed from  the observations of  the function $f$ at the points $x_{i} \in X_N$; see Section~\ref{sec:var_GP}.  If $f$ is indeed a realization of the GP
 with  covariance \eqref{eq:covar}, it follows by the well-known results on
microergodicity \citep[see Chapter 6 in][]{stein1999interpolation}
 that $\widehat{\sigma^2_N } \to \sigma^2$ almost surely  (a.s.)  as $N \to \infty$ and $X_N$ becomes dense in $X$.
Note, however, that such realizations  do not belong to $  H(K)$ a.s. and if $f \in H(K)$ then $\widehat{\sigma^2_N } \to 0$   as $N \to \infty$.
This observation is a consequence of  the useful relation
$
\widehat{\sigma^2_N } = \frac1N \| \mu_N\|_{H(K)}^2
$  \citep[see equation (3.4) in][]{karvonen2020maximum} and the relation $\| \mu_N\|_{H(K)}\leq \| f\|_{H(K)}$, which has already been mentioned.
% and holds for any $f \in H(K)$ and any point set $X_N$.
This is in full agreement with Corollary~\ref{cor:BLUE} below, which  states  that $f\in H(K)$  if and only if $
 \lim_{N \to \infty} N\, \widehat{\sigma^2_N} < \infty
 $ (assuming $X_N$ becomes dense in $X$ as $N \to \infty$). Our numerical studies  in Section~\ref{sec:numer} of this paper  show, however, that
 for finite sample sizes  the asymptotic
 behaviour of $\widehat{\sigma^2_N}$ has much less effect on uncertainty quantification in GP regression than the smoothness of the function $k(x)$
 and its flatness at $x\simeq 0$.
}

The problem of deciding whether  a given function $f$ belongs to the RKHS  $ H(K)$ for a chosen kernel $K$ is well-known in literature, see e.g. Section 3.4 in \cite{paulsen2016introduction}.
%; difficulties in resolving of whether a given function $f$ belongs to $H(K)$  is discussed e.g.  in \cite{steinwart2006explicit}.
This problem is well-studied in the case of kernels with finite number of derivatives,  \citep{ritter2007average,ritter1995multivariate,karvonen2020maximum}. For infinitely differentiable kernels, only the case of Gaussian kernel \eqref{eq:gauss} is well-understood; see \cite{minh2010some} for most advanced results.
 {In a preprint citing the present paper, \cite{karvonen2021non} gives further
new results for analytic translation-invariant kernels with $\lim_{\|x\|\to \infty} k(x)=0$.
 }

\subsection{Sufficient conditions for moment determinacy}
\label{sec:suff}

There  are many
sufficient conditions for moment determinacy of probability measures,
see e.g. \cite{lin2017recent}, \cite{schmudgen2017moment} and \citet[Chapter 11]{stoyanov2013counterexamples}.
The following  sufficient condition for moment determinacy of a measure $\alpha(dt)$ with moments $c_k$
in the Hamburger   moment problem, the so-called {\it Carleman condition}, is one of the most commonly used:
\be
\label{eq:moment_det1} \label{eq:moment_detK}
{\rm (A.1):}\;\;\;\;\;\;\;\;\;\;\;\;\;\;\;\;\;\;\;\;\;\;\;\;\;\;\;\;\;\;\;\;\;\;\;\;\;\;\;\;\;\;\;\;\;\;\;\;\sum_{n=1}^\infty  c_{2n}^{-1/(2n)}= \infty\,. \;\;\;\;\;\;\;\;\;\;\;\;\;\;\;\;\;\;\;\;\;\;\;\;\;\;\;\;\;\;\;\;\;\;\;\;\;\;\;\;\;\;\;\;\; \ee
 {Note that if a  probability measure $\alpha(dt)$ has finite moments $c_{j}$ for all $j=0,1, \ldots$  and satisfies~\eqref{eq:moment_detK},
then for any $m=0,1, \ldots$ the measure $\alpha_m(dt)=t^{2m}\alpha(dt)/c_{2m}$ of Theorem~\ref{th:one_dim1} with moments $\int t^j \alpha_m(dt)= c_{2m+j}/c_{2m}$ also satisfies~\eqref{eq:moment_detK}.
}

 {Less known than  \eqref{eq:moment_detK} is the following  sufficient condition for the moment determinacy (in the  Hamburger   moment problem) of a measure $\alpha(dt)$:
\bea \label{eq:Chow}
{\rm (A.2):}\;\;\;\;\;\;\;\;\;\;\;\;\;\;\;\;\;\;\;\;\;\;\;\; \exists \, \varepsilon>0 \;\mbox{ such that} \;\sum_{n=1}^\infty { |c_{n}| x^{n} \over n! } < \infty\;\;\; \mbox{for all $|x| < \varepsilon$,}\;\;\;\;\;\;\;\;\;\;\;\;\;\;\;\;\;\;
\eea
 see Theorem 30.1 in \cite{billingsley2008probability}. This condition is clearly equivalent to the  assumption (A.3): the
 random variable (r.v.) $\xi_\alpha$ with distribution $\alpha$ has moment generating function and, see Theorem~1 in \cite{lin2017recent}, to (A.4):  $\limsup_{n\to \infty} \frac1{2n} c_{2n}^{1/(2n)}< \infty; $ this is one of the most known and easiest to verify sufficient conditions for moment determinacy.   Theorem 2.13  in \cite{wainwright2019high} yields that the assumptions (A.2)--(A.4) are  also equivalent to the assumption (A.5): the r.v. $\xi_\alpha$ is sub-exponential. In view of this assumption, the tail behaviour of $\alpha(dt)$ is a natural indicator of the degree  of moment-determinacy of $\alpha$.
}

 {
The  conditions (A.2)--(A.5) are stronger than the Carleman condition \eqref{eq:moment_detK} in the sense that  any of them  implies~\eqref{eq:moment_detK}. A technique of constructing  measures $\alpha$ satisfying the Carleman condition \eqref{eq:moment_detK}
but violating $\limsup_{n\to \infty} \frac1{2n} c_{2n}^{1/(2n)}< \infty$, and hence all other assumptions in (A.2)--(A.5), is given in  \cite{stoyanov2013counterexamples}, Section 11.9.
}

In view of \eqref{eq:spectral1} and symmetry of $k(\cdot)$, $|c_{n}|=|k^{(n)}(0)|$ for all integers $n =0,1, \ldots$ This yields that (A.2) coincides with  the assumption that $k(\cdot)$ is a real analytic function at the point $0$, see Definition 1.1.5 and Corollary 1.1.16  in \cite{krantz2002primer}. Summarizing, if $k(\cdot)$ is a real analytic at $0$, then the corresponding spectral measure $\alpha(dt)$ is moment-determinant in the Hamburger sense.  {It is not difficult, however, to construct moment-determinant (in the Hamburger sense) spectral measures $\alpha(dt)$ which do not satisfy (A.2); see Sections 11.9 and 11.10 in \cite{stoyanov2013counterexamples}. For such measures, the corresponding functions $k(\cdot)$ are symmetric, positive definite and infinitely differentiable but not  real analytic at $0$. Nevertheless, for practical purposes the class of functions which are symmetric, positive definite and  real analytic can be considered as the main class of functions $k(\cdot)$, for which the corresponding spectral measures are moment determinant in the Hamburger sense.} Note also that if the function  $k(\cdot)$ is real analytic then all functions in $H(K)$ have to be analytic too; see \cite{sun2008reproducing}.

{Let us give four examples of kernels $K(x,y)=\sigma^2 k(x-y)$ whose spectral measures $\alpha(dt) $ satisfy the Carleman condition \eqref{eq:moment_detK} and hence assumptions of Theorems~\ref{th:one_dim},~\ref{th:one_dim1}  and Corollary \ref{cor:1} ({we will return to these kernels    in Sections~\ref{sec4} and \ref{sec:numer}}). The  spectral density  corresponding to the spectral   measure $\alpha(\cdot) $ will be denoted by $\varphi(\cdot)$, so that $\alpha(dt)=\varphi(t)dt$.}
 { In  examples (E.1)--(E.4) below, we also provide asymptotic expressions for $C_{N}(\alpha)=\sum_{n=1}^N  c_{2n}^{-1/(2n)}$. In view of
 \eqref{eq:moment_det1},
 the rate of divergence of $C_{N}(\alpha)$ (as $N\to \infty$) can also be used to characterize  the degree of moment-determinacy of $\alpha$ (additionally to the tail behaviour of $\alpha$). }

\begin{itemize}
  \item[(E.1)]
  {Gaussian kernel: $k(x)=\exp\{-\lambda x^2/2\}$,   $\varphi_{{\lambda}}(t)=\frac1{\sqrt{2\pi \lambda}}\exp\{-t^2/(2\lambda)\}$ ($t \in \mathbb{R}$),
  $c_{2n}= \lambda^n (2n-1)!!$, $ C_{N}(\alpha)= \sum_{n=1}^{N} c_{2n}^{-1/(2n)} = \sqrt{N/\lambda }(1+ o(1)) ,\, N \to \infty$.}

  \item[(E.2)]
  {Cauchy kernel:
  $k(x)=1/(1+x^2/\lambda^2)$, $\varphi_{{\lambda}}(t)=\frac{ \lambda}{2}\exp\{-{\lambda}|t|\}$  ($t \in \mathbb{R}$),
  $c_{2n}= \lambda^{-n} (2n)!$, $  c_{2n}^{-1/(2n)} =  e\sqrt{\lambda} /n +o(1/n) $, so that
   $ C_{N}(\alpha)  = e\sqrt{\lambda} \log N (1+ o(1)) $ as
   $N \to \infty$.
  }
  \item[(E.3)]
  {
  The kernel whose spectral density is a symmetric Beta-density: \\${\varphi}_a(t)={2^{-2 a -1} } (1-t^2)^{a}/{B(a+1 , a+1 )}$ with $a>-1$ and  $t \in [-1,1]$. \\Here  we have $ C_{N}(\alpha)  = {\rm const} \cdot N (1+ o(1)) $ as
   $N \to \infty$.
  %Two special cases of kernels  in this family are given at the end of  Section~\ref{sec:beta_sp}.
}
 \item[(E.4)]  {$k(x)=\cos(\lambda x)$ with $\lambda \neq 0$; here the spectral measure $\alpha$ is concentrated at 
 $\pm \lambda$ with masses $1/2$ yielding $c_{2n}= \lambda^{2n} $, $c_{2n}^{-1/(2n)}=1/\lambda$ and $ C_{N}(\alpha)=N/\lambda$  for all $N$.}
\end{itemize}

%
% A  well-known sufficient condition for    indeterminacy in   Hamburger   moment problem is the so-called {\it Krein condition}
%\be
%\label{eq:moment_det}
%\int_{-\infty}^{\infty} \frac{-\log \varphi (t)}{1+t^2}< \infty\,
%\ee
% applicable for absolutely continuous measures {$\alpha(dt)=\varphi (t) dt$}, for  which   all moments exist.

%A very clear exposition of the topic of moment  determinacy along with a rich collection of examples of moment-determinant and moment-indeterminant measures is contained in  Chapter~11 of the monograph of \cite{stoyanov2013counterexamples}.

%Two important examples of kernels satisfying conditions of Theorems~\ref{th:one_dim} and~\ref{th:one_dim1} are the Gaussian   kernel
% defined by \eqref{eq:gauss} and having the  spectral density
%  $$\varphi(t)=\frac1{2\sqrt{\pi \lambda}}\exp\{-t^2/(4\lambda)\},  \; t \in \mathbb{R} ,$$
%  and the Cauchy kernel defined by  $k(x)=1/(1+x^2/\lambda^2)$ with $\lambda>0$; the corresponding     spectral density is
%   $$\varphi(t)=\frac{ \lambda}{2}\exp\{-{\lambda}|t|\}\, , \; t \in \mathbb{R} .$$
%   In both cases, one can use   \eqref{eq:moment_detK} for checking that the moment problem is determinate.
%Other examples of kernels satisfying conditions of Theorems~\ref{th:one_dim} and~\ref{th:one_dim1} can be constructed using
%moment-determinate spectral measures; see also  Section~\ref{sec:beta_sp}.
%
%

\subsection{Main steps in  the proofs
%proving  Theorems~\ref{th:one_dim},\ref{th:one_dim1}
and the structure of the remaining part  of the paper}
\label{sec:s}

{
Section \ref{sec:2} is devoted to proving Theorems  \ref{th:one_dim} and \ref{th:one_dim1}; the proofs are given in several steps. The main idea in our approach is to relate
the problem of interest to properties of the best linear unbiased estimate (BLUE)  in linear regression models, which will be worked out in Sections \ref{sec:blue} and \ref{sec:blue1}.
Sections \ref{sec:m_determ}  and \ref{sec:relation1}  provide
different characterizations of the moment determinacy of   spectral measures and finally the proofs
will be completed in Section \ref{sec:proof1}. %and \ref{sec:proof2}.
We now explain the different steps in more detail.
}

{In Section \ref{sec:blue} we consider  a one-parameter linear regression model $y(x)\!=\!\theta f(x)\! +\! \varepsilon(x)$ with $\mathbb{E}\varepsilon(x)\varepsilon(x')\!=\!K(x,x')$ and a regression function $f\!\in \! H(K)$ and show that in this case $\hat\theta_{BLUE} $, the BLUE of $\theta$, exists and its  variance  is strictly positive, see Lemma~\ref{lem:1a}.} We also show that  in the case $f \! \notin \! H(K)$, the BLUE does not exist and
   establish in Lemma~\ref{lem:1} that  for proving $f \! \notin \! H(K)$, it is sufficient to construct a sequence of linear unbiased estimators $\hat\theta _{n}$ of the unknown parameter with variances tending to $0.$
{Such a sequence is constructed in Section  \ref{sec:blue1} for the  location scale model and an explicit expression
for the variance of these estimators in terms of the ratio of determinants
\be
\label{ratio}
{\rm var}(\hat\theta_n) = \frac{{\rm det}  (  c_{2(i+j)} )_{i,j=0}^n}{ {\rm det}(  c_{2(i+j)} )_{i,j=1}^n }\,
\ee
of Hankel-type  matrices of the moments of the spectral measure is derived in Lemma~\ref{lem:10}.
In Section \ref{sec:m_determ} we establish several properties of moment-determinant symmetric measures which we use in
Section  \ref{sec:relation1} for building up an equivalence between the moment determinacy of the spectral measures and the statement that the sequence \eqref{ratio} converges  to zero.
  This is arguably the most important step in the proof of both theorems (see  Lemma~\ref{lem:ham}).
Finally,   these results are combined in the proofs of     Theorem \ref{th:one_dim} and \ref{th:one_dim1}
 in Section \ref{sec:proof1}. } %  and  \ref{sec:proof2}, respectively. }

{
In Section~\ref{sec:3} we consider several extensions  and  interpretations
of the main results.
In Section~\ref{sec:finite_support} we consider   spectral measures with finite support,
%\HD{{\bf delete? } In Section~\ref{sec:finite_mass} we demonstrate that in the case when the spectral measure has positive mass at the point 0 the proof of Theorem~\ref{th:one_dim} fails in the sense that the variances of the constructed estimators $\hat\theta _{n}$ do not converge to zero.}
while  Section~\ref{sec:multivar}   discusses  the multivariate case. This discussion is continued  in Sections~\ref{discret}--\ref{sec:var_GP},
where we also consider  general metric spaces.
In Section~\ref{discret} we explain a   technique of characterizing the inclusion $f \in H(K)$ via suitable discretization of the set $X$ and show that
$1/\|f\|_{H(K)}$ is the limit of variances of the related discrete BLUEs. These results are used in Section~\ref{sec:mass0}, where we prove that the constant function belongs to $H(K)$ if and only if
the spectral measure has positive  mass at $0$.
In Section~\ref{sec:var_GP} we show that   the problem of parameter estimation in a one-parameter regression model is equivalent to the problem of estimating      the variance of a Gaussian process (field). Thus we are able to relate our findings  to  the  estimation problems considered in \cite{xu2017maximum} and \cite{karvonen2020maximum}.
 In Section~\ref{sec:3.1} we return to the one-dimensional case and give an interpretation of  Theorem~\ref{th:one_dim}
in terms of the   $L_2$-error of the best approximation of  a constant function by  polynomials of the form $a_1t^2+a_1t^4+ \ldots +  a_nt^{2n}$.
}

In Section~\ref{sec4}, for two specific classes of kernels we derive explicit results on the rates of convergence to $0$ of the  ratio of determinants \eqref{ratio}. In the case of Gaussian  kernel \eqref{eq:gauss}, we detail and improve one of the  asymptotic expansions
of Theorem~3.3 in \cite{xu2017maximum}. {Finally, in Section~\ref{sec:numer}  we discuss results of a  numerical study for uncertainty quantification in GP regression in  relation to the theoretical results of this paper. }

\section{Parameter estimation, moment determinacy and proofs of main results}
\label{sec:2}
\def\theequation{2.\arabic{equation}}
\setcounter{equation}{0}
 %Theorems \ref{th:one_dim} and \ref{th:one_dim1} will be proved in Sections \ref{sec:proof1} and \ref{sec:proof2},    correspondingly. In these proofs,
%we
%shall use several properties of statistical estimates in regression models, which we will introduce and discuss  in Sections \ref{sec:blue} and \ref{sec:blue1}. We shall also need different characterizations of the moment determinacy of the spectral measure; this will be considered in Sections \ref{sec:m_determ} and \ref{sec:relation1}.

\subsection{{BLUE in a one-parameter regression model %with $f\!\in \! H(K)$ and the case $f \! \notin \! H(K)$
 }}
\label{sec:blue}

\medskip

Consider a one-parameter regression model  with stationary correlated errors:
\be \label{eq:gen_model01}
\mbox{$y(x)\!=\!\theta f(x)\! +\! \varepsilon(x)$, $x\!\in \!X $, $\mathbb{E}\varepsilon(x)\! =\!0,$ $\,\mathbb{E}\varepsilon(x)\varepsilon(x')\!=\!k(x\!-\!x')$.}
\ee
Here $\theta$ is a  scalar parameter, $f\!: \! X\! \to \!\Real $ is a given regression  function and  $k(\cdot)$ is an infinitely differentiable positive definite function
with $k(0)=1$  making the kernel $K(\cdot,\cdot)$  defined by $K(x,y)=k(k-y)$ an infinitely differentiable correlation kernel.
For constructing  estimators of the parameter $\theta$,  the observations of the process $\{y(x)|x \in X\}$ along with observations of all of its derivatives
$\{y^{(k)}(x)|x \in X\}$, $k=1,2, \ldots$, can be used.

{An estimator $ \hat \theta$ for the parameter $\theta$ is called linear, if it is a linear function of the observations (in our case of the process and its derivatives).
An unbiased estimator satisfies $\mathbb{E} [\hat \theta ] = \theta $ for all $\theta$.
The  best linear unbiased estimator (BLUE) of $\theta$ is defined as an unbiased estimator $\hat\theta_{BLUE} $ such that $\mbox{var}(\hat\theta_{BLUE}) \leq \mbox{var}(\hat\theta)$, where $\hat\theta$ is any  linear unbiased estimator  of $\theta$.}
If the kernel $K$ is differentiable and the BLUE exists, then for its computation all available derivatives of $y(x)$ are used, see \cite{dette2019blue}. In general, the BLUE may not exist but the next lemma shows that it does exist when $f \in H(K)$.
\medskip

\begin{lemma}
\label{lem:1a}
If  $f \in H(K)$, then the  BLUE $\hat\theta_{BLUE} $ in  model \eqref{eq:gen_model01}   exists and
$$
\mbox{\rm var}(\hat\theta_{BLUE})= 1/
 \| f \|_{H(K)}>0.
 $$
\end{lemma}

\medskip
The statement of lemma follows from Theorem 6C (p. 975) of \cite{parzen1961approach}. Formally, only the case  $X=[0,1]$ is considered in \cite{parzen1961approach}, but Parzen's proof does not use the structure of $X$ and is  therefore valid
 for  a general metric space $X$.
%  \hfill $\Box$\\

\medskip

\begin{lemma}
\label{lem:1} If   there exists a sequence    of    linear unbiased estimators $(\hat \theta_n)_{n  \in \mathbb{N}}$   of $\theta$ in model \eqref{eq:gen_model01}
such that ${\rm var}(\hat\theta_n)\to 0 $ as $n \to \infty$, then
 $f \notin H(K)$.
\end{lemma}

\medskip

{\bf Proof.} Assume that $f \in H(K)$. By Lemma~\ref{lem:1a}, the continuous BLUE $\hat\theta_{BLUE} $   exists and var$(\hat\theta_{BLUE})= 1/ \| f \|_{H(K)}>0$. From the definition of the BLUE,
 ${\rm var}(\hat\theta_{n}) \geq  {\rm var}(\hat\theta_{BLUE})>0$
  for all $n \in \mathbb{N}$. We have arrived at a contradiction and hence $f \notin H(K)$. \hfill $\Box$\\

%
%
%\subsection{A one-parameter regression model with $f \notin H(K)$ }
%Consider a one-parameter regression model  and correlated errors:
%\be \label{eq:gen_model0}
%\mbox{$y(x)=\theta f(x) + \varepsilon(x)$, $t\in X $, $\mathbb{E}\varepsilon(x) =0,$ $\;\mathbb{E}\varepsilon(x)\varepsilon(x')=K(x,t')$.}
%\ee
%Here $\theta$ is a  scalar parameter, $f(\cdot)$ is a given function on $X$ and $K(x,t')$ is a correlation  kernel.
%
%\begin{lemma}
%\label{lem:1}
%Assume  there exits a family $\hat\theta_n$ (indexed by $n$) of  linear unbiased estimators of $\theta$ such  ${\rm var}(\hat\theta_n)\to 0 $ as $n \to \infty$. Then
% $f \notin H(K)$.
%\end{lemma} \footnote{Parzen has it for $X=[0,1]$, us too, but I think it is trivial to extend it to general $X$}
%
%{\bf Proof.} Assume $f \in H(K)$. In this case, the continuous BLUE $\hat\theta_{BLUE} $ in the model \eqref{eq:gen_model0}   exists and var$(\hat\theta_{BLUE})= 1/ \| f \|_{H(K)}>0$, see \cite{dette2019blue} and \cite[Sect. 6]{parzen1961approach}.
%
%
%
%Let $\hat\theta_{n,BLUE}$ be the BLUE constructed on the base of the same observations as $\hat\theta_n$. Then
%var$(\hat\theta_{n,BLUE})\leq {\rm var}(\hat\theta_n) $ and therefore
% var$(\hat\theta_{n,BLUE}) \to 0$ as $n \to \infty$.
%On the other hand,  var$(\hat\theta_{n,BLUE}) \geq $ var$(\hat\theta_{BLUE})$ for all $n$. We have arrived at a contradiction and hence $f \notin H(K)$. \hfill $\Box$\\
%
\subsection{A family of estimators $\hat\theta_{n}$ in the location scale model}
\label{sec:blue1}

Consider the location scale model
 \be
 \label{eq:loc-sc}
 y(x)\!=\!\theta\!+\!\ve(x)\, , \;x \! \in \! X \subset \! \mathbb{R},\;\; \mathbb{E}\varepsilon(x) \! =\!0,\;\mathbb{E}\varepsilon(x)\varepsilon(x')\!=\!k(x\!-\!x'),
 \ee
 where $k (\cdot)$ is an infinitely differentiable at 0 positive definite function.
 Choose any {interior} point $x_0 \in X$  and set $\ve_0=\ve(x_0)$.
 For construction of the estimator
  $\hat\theta_{n}$, which we will apply in Lemma \ref{lem:1}, we use the following $n+1$ observations: the observation
 $
 y(x_0)=\theta +\ve_0\,
 $ at the point $x_0$
and~$n$ mean-square derivatives of the process  $y$ at the point $x_0$:
 \be
 \label{eq:obs_n0}
 \ve_j =y^{(j)}(x_0) = \frac {d^j y(x)}{d x^j} \bigg|_{x=x_0}  = \frac {d^j \ve(x)}{d x^j} \bigg|_{x=x_0} \, , \;\;\;\; j=2,4,\ldots, 2n\,.
 \ee

{As discussed in Section~\ref{sec:suff}, for the main class of kernels of interest, the RKHS $H(K)$ is a subset of analytic functions. If we observe $y(x)$ everywhere on [0,1], then, since $f \in H(K)$
and $\varepsilon(x)$  are analytic, we also know all $y^{ (k)} (x)$ for all $x \in [0,1]$ and
any integer $k \geq  0$. Again, because of the analyticity, observing $y(x)$ everywhere
on $X = [0,1]$ is the same as observing $y^{ (k)} (x)$ for any $x \in [0,1]$ and any
integer $k \geq 0$. This yields that in practice we do not need to directly observe  $y^{ (k)} (\cdot)$ for constructing estimators of $\theta$.}

The following result provides a necessary and sufficient condition for the existence of the derivatives.
For a proof, see  page 164 (Section 12) in  \cite{yaglom1987correlation}.

\medskip

\begin{lemma}
\label{lemma5}
{
Let $x_0$ be an interior point of $X$. The mean-square derivative
$\ve_j= {d^j \ve(x)}/{d x^j} \big|_{x=x_0}$ of the stationary process
 $\{\varepsilon(x)| x \in X \}$  in \eqref{eq:loc-sc}  at the point $x_{0}$
 exists if and only if $c_{2j}< \infty$, where}
\be
\label{eq:momemts}
c_{2j}= \int_{-\infty}^{\infty} t^{2j} \alpha(dt)=(-1)^j \frac { \partial^{2j}} {\partial u^{2j} } k(u)\bigg|_{u=0}
\ee
 is the $2j$-th moment of the spectral measure $\alpha$ corresponding to the kernel $k$ in Bochner's theorem.
\end{lemma}
\medskip

As we have assumed that the kernel $k (\cdot)$ is infinitely differentiable at  0, all moments $c_j$ ($j=0,1, \ldots$) exist. As an immediate consequence of the existence of all moments and the representation \eqref{eq:spectral},
 for the random variables $\varepsilon_j$ defined in \eqref{eq:obs_n0}, we obtain by Lemma~\ref{lemma5}
%\begin{lemma} \label{lemma5}
\be
\label{eq:deriv}
\mathbb{E}\ve_i \ve_j = \frac { \partial^{i+j}} {\partial x^{i} y^{j}} k(x-y) \bigg|_{x,y=x_0} = (-1)^{i+j}
\frac { \partial^{i+j}} {\partial u^{i+j} } k(u)\bigg|_{u=0} = c_{i+j} \, .
\ee
 for all $i,j=0,1, \ldots$
Note that  all  derivatives ${ \partial^{m}} /{\partial u^{m} } k(u)\big|_{u=0}$ of  odd  order $m$ vanish as the function $k (\cdot)$ is symmetric around the point $0$.

{
Next, we introduce the random variables $\delta_i=(-1)^i \ve_{2i}$, $i=0,1, \ldots$. The observations~\eqref{eq:obs_n0}
used for constructing  the discrete BLUE in model \eqref{eq:loc-sc} can then be rewritten as
$$y_0=y(x_0)= \theta + \ve_0,\; y_1=y^{(2)}(x_0) = \delta_1,\; \ldots,\; y_n=y^{(2n)}(x_0)=\delta_n.$$
Moreover,
the covariance matrix of the vector $(\delta_0,\delta_1, \ldots, \delta_n    )^\top$ is  the Hankel matrix
\be \label{eq:Cn}
C_n =\left(\mathbb{E}\delta_i \delta_j\right)_{i,j=0}^n=( c_{2(i+j)})_{i,j=0}^n\,,
\ee where  $c_2, \ldots, c_{2n}$ are the moments defined in \eqref{eq:momemts}.
}

%This yields, in particular, that for all odd $j$, $E \ve_0 \ve_j=0$ and therefore such $\ve_j$ cannot be used for predicting $\ve_0$ and therefore for constructing the estimator $\hat\theta_n$.
%The estimators $\hat\theta_n$  are constructed as the discrete BLUEs
%of $\theta$ in the model \eqref{eq:loc-sc} using the  observations
%$$y_0=y(x_0)= \theta + \ve_0,\; y_1=y^{(2)}(x_0) = \delta_1,\; \ldots,\; y_n=y^{(2n)}(x_0)=\delta_n.$$

{Assume that the spectral measure $\alpha(dt)$ has infinite support. In this case,   the matrices $C_n$ are positive definite for all  $n=0,1, \ldots$}
 \citep[see, for example,   Proposition 3.11 in ][]{schmudgen2017moment} and  the discrete BLUE is obtained as
\be \label{eq:estim}
 \hat\theta_{n}= \frac{ e_{0,n} ^\top C_n^{-1} Y_n}{ e_{0,n} ^\top C_n^{-1} e_{0,n}},
\ee
where $Y_n=(y_0,y_1,\ldots, y_n)^\top $ and $e_{0,n}=(1,0,\ldots, 0)^\top \in \mathbb{R}^{n+1}$ denotes the first coordinate vector in $\mathbb{R}^{n+1}$.

\medskip

\begin{lemma}
\label{lem:10}
The variance of the estimator  \eqref{eq:estim}
is
\be
\label{eq:var}
{\rm var}(\hat\theta_n)=\frac1{e_{0,n}^\top C_n^{-1} e_{0,n} }= \frac{H_n}{G_n}\, ,
\ee
where $H_n$ and $G_n$ are the determinants
\be\label{eq:H_n}
H_n= {\rm det} \left( C_n  \right)= {\rm det} \left[ \left(  c_{2(i+j)} \right)_{i,j=0}^n \right] \, , \;\;\; G_n= {\rm det} \left[ \left(  c_{2(i+j)} \right)_{i,j=1}^n \right]\, .
\ee
\end{lemma}

\medskip

{\bf Proof.} The expression \eqref{eq:var} follows from the standard formula
${\rm var}( \hat\theta_{n})=1/(e_{0,n} ^\top C_n^{-1} e_{0,n})$ for the
variance of the BLUE and   Cram\'{e}r's rule for computing elements of a matrix inverse; in our case, $e_{0,n} ^\top C_n^{-1} e_{0,n}$ coincides with the top-left element of the matrix $C_n^{-1}$.
\hfill $\Box$\\

{Observing Lemma \ref{lem:1} we conclude that a non-vanishing constant function does not belong to $H(K)$ if
$
\lim_{n \to \infty} {H_n}/{G_n}=0.
$
In the following sections we relate this condition to the moment determinacy of the spectral measure.}

{
\begin{remark}
{\rm
Let us briefly consider
the case where the spectral measure has a positive mass at the point $0$. Consider   the location scale model \eqref{eq:loc-sc}
and let
\be \label{eq:sm_zero} {\alpha}_\gamma(dt)= (1- \gamma ) {\alpha}(dt)+ \gamma \delta_0(dt)\,  \ee
denote the
spectral measure corresponding to a nonnegative definite and symmetric kernel $k_\gamma$,
 where $0<\gamma <1$, $\delta_0$ is the Dirac measure at the point $0$  and $\alpha(dt)$ is a symmetric probability measure on $\mathbb{R}$ with no mass at 0.
The measure ${\alpha}_\gamma(dt)$ is symmetric around the point $0$ with even moments $\tilde{c}_{0}%=\tilde{b}_0
=1$ and
$$
\tilde{c}_{2j}
%=\tilde{b}_j
= (1- \gamma )  c_{2j}
%b_j
\, , \;\; j=1,2, \ldots
$$
Recall the definition of the matrix $C_n$ in \eqref{eq:Cn} and
define the matrices $$\tilde{C}_n =( \tilde{c}_{2(i+j)})_{i,j=0}^n=\gamma  e_{0,n}e_{0,n}^\top + (1- \gamma )  C_n $$   and the corresponding determinants
\bea
\tilde{H}_n= {\rm det} \, \tilde{C}_n\, , \;\;\;
\tilde{G}_n= {\rm det} \left[ \left(  \tilde c_{2(i+j)} \right)_{i,j=1}^n \right]= (1- \gamma ) ^n G_n\, ,
\eea
where $G_n$ is defined in \eqref{eq:H_n}.
Using standard formulas of linear algebra we obtain
\bea
\tilde{H}_n= {\rm det} \,[\gamma  e_{0,n}e_{0,n}^\top + (1- \gamma )  C_n]=(1- \gamma ) ^{n} \left[(1- \gamma )  + \gamma  e_{0,n}^\top  C_n^{-1}   e_{0,n} \right]\,  H_n \, .
\eea
In accordance with \eqref{eq:var}, the variance of $\tilde\theta_n$, the BLUE of $\theta$ constructed similarly to $\hat\theta_n$ but for the spectral measure ${\alpha}_\gamma(dt)$,
is given by
\bea
\label{eq:var5}
{\rm var}(\tilde\theta_n) &=& \frac1{e_{0,n}^\top \tilde{C}_n^{-1} e_{0,n} } = \frac{\tilde H_n}{\tilde G_n} = \frac{ H_n}{ G_n} \left[ (1- \gamma ) +\gamma e_{0,n}^\top  C_n^{-1}   e_{0,n} \right]\\ &=&{\rm var}(\hat\theta_n) [(1-\gamma ) +\gamma / {\rm var}(\hat\theta_n) ]
= (1- \gamma ) {\rm var}(\hat\theta_n) +\gamma    >0
  \, .
\eea
This  implies that ${\rm var}(\tilde\theta_n)$ cannot  converge to $0$  and  Lemma \ref{lem:1} is not applicable if the spectral measure has a positive mass at the point $0$.\\
In Theorem~\ref{th:3} of Section \ref{sec:mass0}  we will  prove  that for any  compact set $X \subset \mathbb{R}^d$ the constant functions indeed  belong to $H(K)$,
if the spectral measure has a positive mass at the point~$0$.
}
\end{remark}
}

\subsection{Moment-determinacy of the spectral measure}
\label{sec:m_determ}

\medskip

Consider the spectral measure $\alpha$  introduced in equation \eqref{eq:spectral}. As a spectral measure, $\alpha$ is a symmetric measure (around $0$) on the real line
and we have assumed that $\alpha$ does not have a positive mass at the point $0$.
%; this assumption is equivalent to $k(x)\to 0$ as $|x|\to \infty$. \HD{Please give a reference!}
Moreover, we have assumed $k(0)=1$ making $\alpha$  a probability distribution.
{In the following  we relate $\alpha$  to a (unique) measure on the   nonnegative axis $[0,\infty)$.
Loosely speaking, if a real valued random variable  $\xi$ has distribution $\alpha(dt)$, then $\alpha_+(dt)$ is the distribution of the random variable
 $\xi^2$. In the opposite direction, if the nonnegative random variable $\eta$ has distribution
$\alpha_+(dt)$, then $\pm \sqrt{\eta}$ has distribution
$\alpha(dt)$, where $\pm$ denotes a random sign.}

{For a more formal   construction  we follow the arguments in Section 3.3 of \cite{schmudgen2017moment} and denote by
 ${\cal B}$   the  Borel sigma field on $\mathbb{R}$, define
 $\tau:\mathbb{R} \to [0, \infty); \tau(x) =x^2$ and $\kappa:[0, \infty) \to \mathbb{R}$, $\kappa(x) = \sqrt{x}$. Then for any symmetric (Radon) measure $\alpha$ on $ {\cal B} $, the measure $\alpha_{+}$  defined by
\begin{equation}\label{h1}
  \alpha_{+}(B ) = \alpha(\tau^{-1}(B )) \qquad B  \in {\cal B} \cap [0,\infty)
\end{equation}
defines a measure on $ {\cal B}  \cap [0, \infty )$. Conversely, if $\alpha_{+}$ is a measure on ${\cal B} \cap  [0,\infty)$, then
\begin{equation}\label{h2}
  \alpha(B ) = \frac {1}{2} \big(\alpha_{+} (\kappa^{-1}(B )) + \alpha_{+}((-\kappa ) ^{-1}(B ))\big)
\end{equation}
defines a symmetric measure on ${\cal B} $. It now follows from Theorem 3.17 in \cite{schmudgen2017moment}  that
the relations \eqref{h1} and \eqref{h2} define a bijection from the set of all symmetric measures on $\mathbb{R}$ onto the set of all measures on $[0,\infty)$.}

The even moments of  a symmetric probability measure  $\alpha$  on $ {\cal B} $
are related to the moments of the measure   $\alpha_{+}$ from \eqref{h1} by
\be
\label{newmoments}
c_{2j}=   \int_{-\infty}^\infty t^{2j}\alpha(dt) =2 \int_{0}^\infty t^{2j}\alpha(dt)
= \int_0^\infty t^{j} d \alpha_{+} (t)=   b_j , \;\;   j\in \mathbb{N}\; ,
\ee
and as a consequence  the determinants $H_n$ and $G_n$ in \eqref{eq:H_n} can be  represented as
\be\label{eq:H_n0}
H_n= {\rm det} \left[ \left(  b_{i+j} \right)_{i,j=0}^n \right]\, , \;\;\; G_n= {\rm det} \left[ \left(  b_{i+j} \right)_{i,j=1}^n \right]\, .
\ee

{Similarly to the case of the Hamburger  moment problem, the {\it Stieltjes moment problem}
is to give necessary and sufficient conditions such that
 a real  sequence $(b_{j})_{j\in \N}$ is in fact a sequence of moments of a  measure $\alpha_{+}(dt)$ on the Borel sets of $[0,\infty)$; that
is  $ b_j= \int_0^\infty t^{j} d \alpha_{+} (t)$ for all $j  \in \N_{0} $.   The Stieltjes moment  problem is determinate
  if the  sequence  of moments  $(b_{j})_{j\in \N}$  %$\{ \tilde{c}_0,\tilde{c}_1, \ldots\}$
determines the measure $\alpha_{+}(dt)$ uniquely. For a proof   of the following result, which relates the Hamburger and Stieltjes moment problem,
 see  \citet[Lemma 1]{heyde1963some},     \citet[Proposition 3.19]{schmudgen2017moment} and    \citet[Sect. 11.10]{stoyanov2013counterexamples}.
}

\medskip

\begin{lemma}
\label{lem:equiv} Let   $\alpha$ be a symmetric probability  measure  on $ {\cal B} $.
The Hamburger moment problem for  $\alpha$ is determinate if and only if the Stieltjes moment problems for the measure
$\alpha_{+}$ defined by \eqref{h1} is determinate.
\end{lemma}

\medskip

Note that  for the  equivalence in Lemma~\ref{lem:equiv} to hold, the assumption  that $\alpha$ does not have mass at $0$ is not required. This assumption, however, is needed in the next lemma.

\medskip

\begin{lemma}
\label{lem:equiv2} Let   $\alpha$ be a symmetric probability  measure  on $ {\cal B} $ with no mass at the point~$0$.
The Hamburger moment problem for  $\alpha$ is determinate if and only if the Hamburger  moment problem for the measure
$\alpha_{+}$ defined by \eqref{h1} is determinate.
\end{lemma}

\medskip

{\bf Proof.} Using the result of Theorem A in \cite{heyde1963some} (see also \cite[p.113]{stoyanov2013counterexamples} and \cite[Remark 2.12]{schmudgen2017moment}), if the Stieltjes  moment problems for the measure
$\alpha_{+}$ is determinate and the measure $\alpha_{+}$ has no mass at~0, then the Hamburger  moment problems for this measure
is also  determinate. From Lemma~\ref{lem:equiv}, the required  equivalence follows.
\hfill $\Box$

\subsection{Relating  moment-determinacy of the measure $\alpha_{+}$ to ${\rm var}(\hat\theta_n)$}
\label{sec:relation1}

\begin{lemma} \label{lem:ham} Let   $\alpha$ be a symmetric probability  measure  on $ {\cal B} $ with infinite support and no mass at
the point $0$.
The Hamburger moment problem for the measure $\alpha_{+}$ defined by \eqref{h1} is determinate if and only
if $H_n/G_n \to 0$ as $n \to \infty$,
where the determinants $H_n$ and $G_n $ are  defined in \eqref{eq:H_n0}.
\end{lemma}

\medskip

{\bf Proof.}
(i) Assume that the  moment problem for the measure $\alpha_{+}$  is determinate.
Let $\mathcal{P}_n$ denote the class of all polynomials of degree $n$ and define
$$
\rho_n(t_0) = \min \left \{ \int_{\mathbb{R}} |P_n(t)|^2 \alpha_{+} (dt) \mid P_n \in \mathcal{P}_n, P_n(t_0) =1 \right \}
$$
for any $t_0 \in \mathbb{R}$, which is not a root of the $n$th  orthogonal polynomial with respect to the measure $\alpha_{+}$
\citep[see equation  (2.26) in  Lemma 2.11  of][]{shohat1943problem}.
%\footnote{ Holger:"In the proof of Lemma, end of first paragraph: I think the statement regarding the limit $\rho (t)$ is only correct for those t which are no zeros  of the orthogonal polynomials. This is not important because we only need it for t=0 and zero is definitively not a root of any orthogonal polynomial with respect to G (as G is supported on the half-line $(0, \infty))$." \\
%A: we can modify the sentence but in Th. 2.6 in Sc and Tamarkin $\rho(t)=\lim_{n \to \infty}\rho_n(t)$, $\forall t>0$, I think.
%\HD{\bf I agree, but I would like to carefully defined $\rho_n$ and for this we have to use either two definitions or make an assumption on $x_0$!}
%}
Then
$$\lim_{n \to \infty}\rho_n(t_0)=: \rho(t_0)$$
exists, by  Theorem 2.6  in \cite{shohat1943problem}.
As the point $0$ is not a support point of the measure $\alpha_{+}$ and all roots of the orthogonal polynomials
with respect to the measure $\alpha_{+}$ are located in supp$(\alpha_{+})\subset (0, \infty)$ we have
from Corollary 2.6 in \cite{shohat1943problem}
 that
 $$
 \rho(0)=\lim_{n \to \infty} \rho_n(0)=0.
 $$
Moreover, by the discussion on  p. 72 (middle of the page) in \cite{shohat1943problem}  it follows that $\rho_n(0)$ is exactly the ratio
$H_n/G_n $, where $H_n$ and $G_n $ are the determinants in \eqref{eq:H_n0}.
 Hence the moment determinacy  for the measure $\alpha_{+}$
implies
$H_n/G_n \to 0$ as $n \to \infty$.
\smallskip

(ii) To prove the converse, assume that $H_n/G_n \to 0$ as $n \to \infty$.
Let $\lambda_n$ be the smallest eigenvalue of the matrix $C_n$.
Theorem 1.1 in \cite{berg2002small} states that the condition
$$ \lim_{n \to \infty} \lambda_n = 0$$
 is necessary and sufficient for the moment-determinacy  of the measure $\alpha_{+}$.

From the definition of $\lambda_n$ as the smallest eigenvalue of the matrix $C_n$ and the representation       \eqref{eq:var} it  follows
$$\lambda_n \leq {1 \over  e_{0,n}^\top C_n^{-1} e_{0,n}} =  \frac{H_n}{G_n} = \rho_n(0)$$
for all $n \in \mathbb{N}$  \citep[see also a related discussion in][]{berg2002small}. Therefore, $H_n/G_n \to 0$ as $n \to \infty$ implies $\lambda_n \to 0$ as $n \to \infty$ and this yields the moment determinacy  of the measure $\alpha_{+}$.
\hfill $\Box$\\

\subsection{Proof of Theorem~\ref{th:one_dim}  and  ~\ref{th:one_dim1} %for the case  $f(x)={\rm const}\neq 0$
}
\label{sec:proof1}

\medskip

{\it Proof of Theorem~\ref{th:one_dim}. }
 Use Lemma~\ref{lem:1} with the estimator defined in \eqref{eq:estim}.
By Lemma~\ref{lem:10} the variance of this estimator is given by \eqref{eq:var}.
From Lemma~\ref{lem:ham}, the determinacy of the measure $\alpha_{+}$ is equivalent to ${\rm var}(\hat\theta_n) \to 0$ as $n\to \infty$.
By  Lemma~\ref{lem:equiv}, this is also equivalent to the moment determinacy of the spectral measure~$\alpha$.
\hfill $\Box$\\

\medskip

%\subsection{Proof of Theorem~\ref{th:one_dim1} % for the case when $f(x)$ is a polynomial}
%}
%\label{sec:proof2}
%{\bf Proof.}
{\it Proof of Theorem~\ref{th:one_dim1}. }
Assume that the function $f$ in \eqref{eq:gen_model01} is a  polynomial of degree $ m \!\geq \! 1$. Take
$m$ derivatives of both sides in \eqref{eq:gen_model01}. The model \eqref{eq:gen_model01} thus reduces to
$\label{eq:gen_model0a}
\mbox{$\tilde{y}(x)=\tilde{\theta}  + \tilde{\varepsilon}(x)$, $x\in X $, }
$
where $\tilde{\theta}$ is the new parameter, $\tilde{y}(x)= {y}^{(m)}(x)$ are new observations and
$\tilde{\varepsilon}= {\varepsilon}^{(m)}$ is the new error process.
From  (2.178) in \cite{yaglom1987correlation},
the autocovariance function of the process $\{{\varepsilon}^{(m)}(x)|x \in X\}$ is given by
$$\mathbb{E} \varepsilon^{(m)} (x) \varepsilon^{(m)}(x')={k}_m(x-x') \;\; {\rm with} \;\;{k}_m(x)=(-1)^m {k}^{(2m)}(x). $$
From \eqref{eq:spectral}, the spectral measure associated with the kernel ${k}_m(x-x')$ is
$\alpha_m(dt) =t^{2m} \alpha(dt)/c_{2m}$. Hence, the statement for the case  when $f$ is a polynomial of degree $m \geq 1$ is reduced to the case of the constant function proved in Theorem~\ref{th:one_dim}; this theorem is applicable  as the measure $\alpha_m(dt) $ does not have mass at 0 for any  $m\geq 1$.
\hfill $\Box$\\

\section{Extensions of Theorems \ref{th:one_dim} and \ref{th:one_dim1} and further discussion}
\label{sec:3}
\def\theequation{3.\arabic{equation}}
\setcounter{equation}{0}

In this section we discuss several extensions of the results derived in Sections \ref{sec1} and   \ref{sec:2}. In particular, we consider spectral  measures with positive mass at the point  $0$ and  extends the  results to the multivariate case. Moreover, we  briefly indicate a relation of our results to the optimal approximation of a constant function by  polynomials with no intercept.

\subsection{{Spectral measures with finite support} }
\label{sec:finite_support}
If   the spectral measure $\alpha(dt)$ in \eqref{eq:spectral} has finite support, say
${\cal T}=\{\pm t_1, \ldots, \pm t_m\}$ with $m \geq 1$ and $0<t_1<\ldots<t_m$, then the matrices $C_n$ in \eqref{eq:Cn} are invertible for $n \leq m-1$ but
\begin{equation} \label{hol1}
{\rm det}(C_n) = {\rm det}  (  c_{2(i+j)}  )_{i,j=0}^n  =0 ~\mbox{ for } ~ n \geq m.
\end{equation}
{Consequently, observing Lemma \ref{lem:10} we have in this case
$$
\mbox{var}(\hat \theta_n)=0 \quad \quad \mbox{for } n=m,m+1,\ldots
$$
Therefore, by Lemma \ref{lem:1} a non-vanishing constant function does not belong to $H(K)$ if the corresponding spectral measure has finite support. }
\medskip

  The relation  \eqref{hol1} follows, observing the representation
\begin{equation*}
  C_n = 2 \sum_{i=1}^{m} w_i g(t_i) g^\top (t_i) \in \mathbb{R}^{(n\!+1) \times (n\!+\!1)}
\end{equation*}
 where $g(t) = (1, t^2,\ldots,t^{2n})^\top$ and $w_1,\ldots,w_m$ are the masses  of the measure $\alpha$
 at the points $t_1,\ldots,t_m$. As $C_n$ is a sum of rank one matrices, it is singular whenever $n > m-1$. On the other hand, in the case $m=n+1$ we have by the Vandermond determinant formula
 \begin{equation*}
   {\rm det} \ C_n =  \prod^{n+l}_{i=1} (2 w_i) \prod_{1 \leq i < j\leq  n+1} (t^2_i - t^2_j)^2 > 0~,
 \end{equation*}
 which shows that $C_n$ is nonsingular. Finally, if $m \geq n+1$ we have (in the Loewner ordering)
 \begin{equation*}
   C_n \geq 2  \sum_{i=1}^{n+1} w_i g (t_i) g^\top (t_i)
 \end{equation*}
 where the matrix on the right-hand side is positive definite.

%Observing the representation \eqref{eq:eq_Vv}  we obtain from \eqref{hol1}
% that
% $$
%  V(p_n^*) = {\rm var}(\hat\theta_n) = 0 ~,~~n=m, m+1, m+2, \ldots
%  $$
%  for any minimizer $p_n^*$ of the functional $V$ in ${\cal P}_{n-1}$.

\subsection{{Multivariate case  }}
\label{sec:multivar}
\medskip

Consider the location scale model
\eqref{eq:loc-sc} but assume  that $X $ is a subset of $ \mathbb{R}^d$ with non-empty interior.
Extensions of Theorems \ref{th:one_dim} and \ref{th:one_dim1} to the  multivariate case, when $d>1$, essentially  follow from the one-dimensional results because it is sufficient to   use      derivatives of  the process $\{ y(x); x \in X \} $ with respect to  one variable for  construction of estimators $(\hat\theta_n)_{n\in  \N}$   { and subsequent  application of Lemma \ref{lem:1}.}
In the following discussion we consider two cases for the kernel $K(x,x^\prime)$ using the notation  $x=(x_1, \ldots, x_d)^\top$, $x^\prime=(x_1^\prime, \ldots, x_d^\prime)^\top$ and  $t=(t_1, \ldots, t_d)^\top$.
We also  denote by  ${x}_{(i)},x^\prime_{(i)}$ and ${t}_{(i)}\in \mathbb{R}^{d-1}$ the vectors $x,x^\prime$ and $t$ with $i$-th component removed, respectively.
\medskip

{\it Case 1:} Assume that $K$ is a product kernel, that is
\be \label{case1}
K(x,x^\prime)=\prod_{j=1}^{d} K_i(x_j,x_j^\prime)\, ,
\ee
where for all $j=1, \ldots , d$ the kernel $K_j$ (defined on  a subset of $\R^2$) satisfies $K_j(x_j,x_j^\prime)=k_j(x_j-x_j^\prime)$  and   $k_j$ is a  non-constant positive definite  function infinitely differentiable at the point $0$. Denote by $\alpha_j(dt_j)$  the spectral measure for $k_j$ and  define  $\alpha(dt)= \alpha_1(dt_1)\cdots \alpha_d(dt_d)$.
To construct the sequence of estimators $(\hat{\theta}_n)_{n\in \N}$  {for the application of Lemma \ref{lem:1}}, we can use the derivatives with respect to the $i$-th coordinate for any $i$. Therefore, Corollary~\ref{cor:1}
can be  generalized as follows.

\begin{corollary}
\label{cor:1m}
Assume that $X \subset \mathbb{R}^d$ and the kernel $K$ has the form \eqref{case1}.
  Then we have the following:
\begin{itemize}
 \item[(a)] If the measure  $\alpha $ has a positive mass at  the point $ 0$,  then the constant functions
%$f(x)={\rm const}\neq 0$, $\forall x \in X$,
belong to $H(K)$.
  \item[(b)]
If  for at least one $i \in \{ 1, \ldots , d\} $ the  Hamburger  moment problem for the measure $ \alpha_i(\cdot)$ is determinate and the measure  $ \alpha_i $ does not have  a positive mass at
the point~$0,$
  then  any non-vanishing  constant function
does not belong to $H(K)$.
  \item[(c)]  If for at least one $i \in \{ 1, \ldots , d\} $ the  Hamburger  moment problem for the measures $t_i^{2m} \alpha_i(dt_i)/c_{2m}$  is  determinate for all $m=0,1,\ldots$, then $  H(K)$ does not contain non-constant  polynomials on $X$.
\end{itemize}
\end{corollary}

Note that the set   $X$ in Corollary~\ref{cor:1m} does not have to be a product of one-dimensional sets.
Moreover, we also point out that the assumption \eqref{case1} can be generalized to kernels of the form
$$
K(x,x^\prime)= k_i(x_i-x_i^\prime) {K}_{(d-1)}({x}_{(i)},x^\prime_{(i)}),
$$
 where ${K}_{(d-1)}(\cdot,\cdot)$ is a positive definite and suitably differentiable kernel on $\mathbb{R}^{d-1}\times  \mathbb{R}^{d-1}$ and $k_i$ is a non-constant positive definite  function infinitely differentiable at the point $0$. % A version of Corollaries~\ref{cor:1} and~\ref{cor:1m} easily follow.

\medskip

{\it Case 2:} The kernel $K$ satisfies
$$
K(x,x^\prime)=k(x-x^\prime),
$$
where $k$ is a positive definite function on $\R^d$. Consider the spectral measure $\alpha(dt)$  corresponding to  $k$ by Bochner`s theorem, that is
\be \label{eq:spectrald}
k(x)= \int_{\mathbb{R}^{d}} e^{i (t_1 x_1+ \ldots+t_d x_d)} \alpha(dt)~,
\ee
and denote  by
$$
\alpha_i(B)=\int_{\mathbb{\R}^{d}} I_B(t_i)  \alpha(dt ) ~,~~ B \in \mathcal{B}  ,
$$
the $i$th  the marginal distribution of  the measure  $\alpha $ ($i=1, \ldots, d$), where $I_B$ denotes the indicator function of the set $B$.
In this case, we can generalize Corollary~\ref{cor:1}
as follows.

\medskip

\begin{corollary}
\label{cor:2}
If  the spectral measure  $\alpha(dt)$ does not have a positive mass at the point~$0$ and
if for at least one  $i \in \{1, \ldots, d\}$ the  Hamburger  moment problems for the measures proportional to $t^{2m} \alpha_i(dt_i)$  are  determinate for
all $m =0,1, \ldots $, then $  H(K)$ does not contain  non-vanishing  polynomials.
\end{corollary}

\medskip

The case when the   spectral measure has positive mass at the point $0$ is treated similarly in one-dimensional and multi-dimensional cases, see Section \ref{sec:mass0}.

\subsection{{Discretization of space and  the limit of    discrete BLUEs}}
\label{discret}

\medskip

{
In Section \ref{sec:mass0}  below we will  prove that  constant functions    belong to $H(K)$ if the spectral measure has positive mass at the point $0$.
The proof requires an   auxiliary  result  which is of own interest and shows that    in the case $f \in H(K)$
 the variance of the continuous BLUE is the limit of the variances of discrete BLUEs, after a suitable discretization of $X$ has been  performed.
}

\begin{lemma} \label{lem:BLUE}
Let $X$ be a compact in $\mathbb{R}^d$,
 $(x_N )_{N  \in \N} $ be a sequence of distinct points in $X$
 such that $f(x_1) \neq 0$ and   \be \label{eq:sup_p}
 \sup_{x \in X} \min_{1\leq i \leq N } \|x-x_i\| \to 0\;\;{\rm as}\;  N  \to \infty\, .\ee
%Let, for a given $n$, $X_n=\{x_1,\ldots,x_n\} \subset X$ \footnote{double index for $x_{i,n}$ ?}
%be a point set and
Let $\widehat{\theta}_{BLUE,N }$ be the BLUE  of
$\theta$ in   model \eqref{eq:gen_model01} from the  observations of $y(x_1) , \ldots , y(x_N ) $.
Then $f \in H(K)$ if and only if var$(\widehat{\theta}_{BLUE,N }) \to c>0$ as $N  \to \infty$.

Moreover, if $f \in H(K)$,  the continuous BLUE $\hat \theta_{BLUE}$   of $\theta$  in model \eqref{eq:gen_model01} exists and
$$c=1/\|f\|_{H(K)}={\rm var}(\hat \theta_{BLUE}).
$$
\end{lemma}

{\bf Proof.}
%Lemma~\ref{lem:1a} and
{Let $X_N  =\{x_1, \ldots , x_N \}$,  $K_N $ denote the restriction of $K$ on $X_N $, and define $H_N =H( K_N )$ as the RKHS corresponding to  the kernel $K_N $.
By Theorem 6 in Section 1.4.2 of  \cite{berlinet2011reproducing} we have  for $f_N =f\big|_{X_N }$, the restriction of $f$ on $X_N $, that   $f_N  \in H_N   $ and
 $
 \|   f_N  \|_{H_N  } \leq \|
 f_{N +1} \|_{H_{N +1} } \leq  \| f \|_{H(K)}\, .
 $}
Consequently,  the sequence of (var$(\widehat{\theta}_{BLUE,N }))_{N  \in \mathbb{N}} { = (1/  \|   f_N  \|_{H_N  })_{N  \in \mathbb{N}} }$ is monotonously decreasing so that
the limit $c= \lim_{N  \to \infty} {\rm var}(\widehat{\theta}_{BLUE,N })\geq 0$ exists for any  $f$. Moreover, ${\rm var}(\widehat{\theta}_{BLUE,N })\geq c $ for all $N  \in \N $.
If  $f \in H(K)$ we have by  Proposition 3.9 in \cite{paulsen2016introduction}  that
$
\lim_{N  \to \infty} {\rm var}(\widehat{\theta}_{BLUE,N }) =c = 1/\|f\|_{H(K)}\, .
$
Conversely,  if var$(\widehat{\theta}_{BLUE,N }) \to c$ as $N  \to \infty$  for some $c>0$, we can use
 the equivalence between (1) and (2) in Theorem 3.11 of  \cite{paulsen2016introduction} to  deduce that $f \in H(K)$.
\hfill $\Box$\\

Recall that  the explicit expression for the variance of the discrete BLUE $\widehat{\theta}_{BLUE,N }$ of Lemma~\ref{lem:BLUE} is given by
\be \label{var_d_BLUE}
{\rm var}(\widehat{\theta}_{BLUE,N })= 1/ F_N ^\top K_N ^{-1} F_N \, ,
\ee
where
\be \label{eq:XW1}
F_N = (f(x_1), \ldots, f(x_N ))^\top\, , \;\;K_N =(K(x_i,x_j))_{i,j=1}^N \, .
\ee
{Since  the  kernel $K(\cdot,\cdot)$ is assumed to be  strictly positive definite, the  matrix $K_N $ is invertible for all $N =1,2, \ldots$ }

\medskip

\subsection{{Spectral measures  with  positive mass at the point $0$}}
\label{sec:mass0}

\medskip
In this section, we investigate the case, where the spectral measure has a    positive  mass at  the point $0$ in more detail. In particular, we show that in this case the constant functions belong to $H(K)$.
To be precise, assume that the covariance kernel  of the error process has the form
\be \label{eq:kernel} K_\gamma(x,x')=\gamma+(1-\gamma)K(x,x')\, , \ee
where $0\leq \gamma<1$ and  $K(x,x')$ is a strictly positive definite kernel on a compact set $X \subset \mathbb{R}^d$. Note that in the particular case $d=1$ and
$K(x,x')=k(x-x')$ with $k$ having the spectral measure $\alpha(dt)$, we obtain the representation \eqref{eq:sm_zero} for the spectral measure~$\alpha_\gamma$.

\medskip

\begin{theorem}
\label{th:3} Let $X \subset \mathbb{R}^d$ be a compact set and assume  the kernel $K_\gamma$ has the form \eqref{eq:kernel} with $0<\gamma<1$. Then
then the  constant functions
%$f(x)={\rm const}\neq 0$, $\forall x \in X$,
belong to  $H(K_\gamma)$.
\end{theorem}
{\bf Proof.}
Consider the location scale model
\be
 \label{eq:loc-sc1}
 y(x)\!=\!\theta\!+\!\ve(x)\, , \;x \! \in \! X,\;\; \mathbb{E}\varepsilon(x) \! =\!0,\;\mathbb{E}\varepsilon(x)\varepsilon(x')\!=\!K_\gamma(x,x') \, .
 \ee
and let
 $(x_n)_{n \in \N}$ denote  a sequence of distinct points in $X$ such that  \eqref{eq:sup_p} is satisfied.
Let $\hat \theta_{m,\gamma} $ be the BLUE of $\theta$ in the  model \eqref{eq:loc-sc1}, constructed on the observations   of
$y(x_1) ,\ldots, y(x_m).$
Define  $W_{m,\gamma}=(K_\gamma(x_i,x_j)_{i,j=1}^m$,  ${Y}_m= (y(x_1),\ldots, y(x_m))^\top $ and ${\bf 1}_m= (1,\ldots, 1)^\top \in \mathbb{R}^m$. As
the covariance kernel $K(x,x')$ is  strictly positive definite, the matrix $W_{m,\gamma}$ is invertible for all $m\geq 1$, $0\leq \gamma<1$. Therefore, the BLUE  is unique  and given  by
$$
\hat \theta_{m,\gamma} = {\bf 1}_m^\top W_{m,\gamma}^{-1} {Y}_m/ {\bf 1}_m^\top  W_{m,\gamma}^{-1}{\bf 1}_m \, .
$$
Its variance is
$$
{\rm var}(\hat \theta_{m,\gamma}) = 1/{\bf 1}_m^\top W_{m,\gamma}^{-1}{\bf 1}_m\, .
$$
For simplicity of notation, denote $\kappa_{m,\gamma}   = {\bf 1}_m^\top W_{m,\gamma}^{-1}{\bf 1}_m=1/{\rm var}(\hat \theta_{m,\gamma})$. The same arguments as given in the proof of
Lemma~\ref{lem:BLUE} show that   for any $0\! \leq \! \gamma \!<\!1$, the sequence $( \kappa_{m,\gamma} )_{m \in \N} $ is monotonously increasing  with  some limit $c_\gamma\!=\! \lim_{m\to \infty} \kappa_{m,\gamma} \in ( 0, \infty] $. Observing the representation
$$
W_{m,\gamma}=(1-\gamma)W_{m,0}+\gamma {\bf 1}_m {\bf 1}_m^\top
$$
(for all $m=1,2, \ldots$ and $0 < \gamma <1$), we have
\bea
W_{m,\gamma}^{-1}=
\frac1 {1-\gamma}  \Big[ W_{m,0}^{-1} - \frac{\gamma }{1-\gamma +\gamma {\bf 1}_m^\top W_{m,0}^{-1} {\bf 1}_m } W_{m,0}^{-1} {\bf 1}_m
{\bf 1}_m^\top W_{m,0}^{-1} \Big]\, .
\eea
This implies
\bea \label{eq:kappan1}
\kappa_{m,\gamma}=
\frac{\kappa_{m,0}}{1- \gamma } \Big[1 - \frac{\gamma\kappa_{m,0} }{1-\gamma +\gamma\kappa_{m,0}  } \Big]=
\frac{\kappa_{m,0}} {1-\gamma +\gamma\kappa_{m,0}   },
\eea
and therefore it follows that
\be \label{eq:kappan}
{\rm var}(\hat \theta_{m,\gamma})= 1/\kappa_{m,\gamma}= \gamma+ (1-\gamma) {\rm var}(\hat \theta_{m,0}) \, .
\ee

Taking the limit (as $m \to \infty$) in \eqref{eq:kappan} we obtain for all $0\! < \!\gamma \!<\!1$:
\bea
 \lim_{m\to \infty} {\rm var}(\hat \theta_{m,\gamma})= \gamma+(1-\gamma)/c_0 \geq \gamma>0\, .
\eea
Lemma~\ref{lem:BLUE} now yields that
the  constant functions
%$f(x)\!=\!{\rm const}\!\neq \!0$ ($\forall x \!\in \!X$)  %
belong to  $H(K_\gamma)$.
 \hfill $\Box$\\

\subsection{Estimation of  the  variance of a Gaussian random field}
\label{sec:var_GP}
\medskip

Let $X \subset \mathbb{R}^d$  be a compact set,  and let $f$   denote of a Gaussian random process (field) on $X$  with a strictly positive definite covariance kernel
$R(x,x')=\sigma^2 K(x,x')$  on $X\times X$, where the kernel $K(x,x')$ is known  but $\sigma^2$ is unknown.
For estimating $\sigma^2$ we assume that one can observe $f$ at $N$ distinct points $x_1,  \ldots, x_N \in X$.
Then it is easy to deduce (see, for example, p.140 in \cite{xu2017maximum}) that the corresponding log-likelihood function is
\be \label{eq:ll}
LL(\sigma^2)= \frac12 \left[ -N   \log(2\pi)-N  \log  ( \sigma^2) -\log( {\rm det}( K_N ) ) -
\frac{1}{\sigma^2}
F_N ^\top K_N  ^{-1}F_N  \right] , \;\;\;\;\;\;
\ee
where  $F_N $ and $K_N $  are defined by \eqref{eq:XW1}. Moreover,  a simple  calculation shows that the
%Maximizing $LL(\sigma^2)$ of \eqref{eq:ll} with respect to $\sigma^2 $, we obtain the
maximum likelihood estimator (MLE) of $\sigma^2 $ is given by
\be
\label{1s}
\widehat{\sigma^2_N  } = \frac1N  F_N ^T K_N ^{-1} F_N \, .
\ee
Comparing \eqref{1s} with \eqref{var_d_BLUE} we get
\be
\label{eq:connect}
\widehat{\sigma^2_N  }=\frac{1}{N \, {\rm var} (\widehat{\theta}_{BLUE,N })}\, ,
\ee
and  by Lemma \ref{lem:BLUE} we obtain  the following corollary.

\medskip

\begin{corollary} \label{cor:BLUE}
{
Let $X$ be a compact set  in $\mathbb{R}^d$, $K$  a strictly positive definite kernel on $X \times X$ and
$f$  a function on $X$. If
 $x_1,x_2,\ldots$  is a sequence of distinct points in $X$ satisfying   \eqref{eq:sup_p}  and
 $\widehat{\sigma^2_N  }$ is  the MLE of $\sigma^2$ constructed from the observations  $f(x_1), \ldots, f(x_N )$  under the assumption that $f$ is a realization of a GP with zero mean and covariance \eqref{eq:covar}, then  $f \in H(K)$  if and only if
 $
 \lim_{N  \to \infty} N \,\widehat{\sigma^2_N } < \infty.
 $
 }
\end{corollary}
\medskip

\subsection{{Best polynomial approximation}  }
\label{sec:3.1}

\medskip

Let   $L_2(\alpha)$ denote the space of square integrable functions with respect to the measure $\alpha (dt ) $ on the real line and define ${\cal P}_{n-1}$ to be the space of
of polynomials of degree $n-1$. For $p \in \mathcal{P}_{n-1}$ we consider the $L_2(\alpha)$-distance
$$
V(p)=\int_{-\infty}^{-\infty} (1-t^2 p(t^2) )^2 \alpha(dt)
$$
between the constant function $g(t) \equiv 1$ and the even polynomial $t^2 p(t^2) $ of degree $2n$ with no intercept. A well know result in approximation theory \citep[see, for example,][p. 15-16]{achieser1956} shows  that
\be
\label{eq:eq_Vv}
 {\rm min}_{ p_n \in {\cal P}_{n-1} } V(p)   =  {{\rm det}  (  c_{2(i+j)}  )_{i,j=0}^n  \over  {\rm det}  (  c_{2(i+j)}  )_{i,j=1}^n  }
% \int_{\T}(1-t^2 p_n(t^2) )^2 \alpha(dt)
= {\rm var}(\hat\theta_n)\, ,
\ee
where  $c_0, c_2,  c_4, \ldots $ are the (even) moments of the spectral measure $\alpha$ defined in \eqref{eq:momemts} and the last equality is a consequence of  Lemma \ref{lem:10}.

From this representation it follows that
${\rm var}(\hat\theta_n)\to 0$ as $n \to \infty$ if and only if non-zero constant functions can be approximated by polynomials
of the form $\tilde{p}_n(t)=t^2 p_n(t^2)$ with arbitrary small error. Moreover,
for any polynomial $p $ on $(-\infty, \infty)$, we have
\bea
 V(p)
%= \int_{-\infty}^{\infty}(1-t^2 p(t^2) )^2 \alpha(dt)
&=&  \int_{0}^{\infty}(1-t p(t) )^2 \alpha_+(dt)= b_2 \int_{0}^{\infty}(1/t- p(t) )^2 {\alpha}_{2,+}(dt),
\eea
where the measure  $ \alpha_+(dt)$ is defined by \eqref{h1}, $ b_2= \int_0^\infty t^{2} d \alpha_{+} (t)$ and ${\alpha}_{2,+}(dt)= t^2 \tilde{\alpha}_+(dt)/b_2$.
{
From Corollary  2.3.3 in   \cite{akhiezer2020classical},  it therefore follows that  the set of all polynomials
${\cal P}_\infty = \cup_{n=0}^\infty {\cal P}_n$
is dense in the space $L_2( [0, \infty) , \nu )$ if the  measure $\nu$  on $ [0, \infty)$ is  the   (unique) solution of    a determinate  Hamburger moment problem.}
As the  function $f(t)=1/t$ belongs to $L_2( (0, \infty), {\alpha}_{2,+})$ we thus obtain from \eqref{eq:eq_Vv}
 another proof of the fact that if $\alpha(dt)$ has no mass at $0$  and
${\alpha}_{2,+}(dt)$ is moment-determinate in the Hamburger sense then ${\rm var}(\hat\theta_n) \to 0$.
Note that   this is almost equivalent to the `if' statement in the important Lemma~\ref{lem:ham}.

\section{Rates of convergence} % of $H_n/G_n$ to 0}
\label{sec4}

\def\theequation{4.\arabic{equation}}
\setcounter{equation}{0}

In this section, we derive for several specific classes of correlation kernels   explicit results on the rate of convergence of the ratio ${\rm var}(\hat\theta_n)=H_n/G_n$, see \eqref{eq:var},  where $H_n$ and $G_n$ are the determinants defined in \eqref{eq:H_n}.

\subsection{Gaussian  kernel}
We first consider the case of the Gaussian  kernel $K(x,x')= \exp\{- \lambda (x-x')^2\}$ with $X \subset \mathbb{R}$ and $\lambda>0$. Assuming for simplicity    $\lambda=1/4$, we obtain that
 the spectral measure  is absolute continuous with density
$$
\varphi(t)=\frac{1}{\sqrt{\pi}} e^{-t^2}, \;\; -\infty<t<\infty\, .
$$
The moments of even order of the measure $\alpha$ are given by
$$
c_{2j}=\int_{-\infty}^\infty t^{2j}\varphi(t) dx=     \int_0^\infty t^{j} g(t) dt
=b_j = 2^j(2j-1)!! \quad j=0,1,\ldots,
$$
where
$
g(y)\!=\!\frac{1}{\sqrt{\pi}} y^{-1/2} e^{-y}, \;y>0.
$
% That is, $b_j$ is $j$-th moment of the distribution with density $g(y)$; $j=0,1,\ldots$.
{From \eqref{eq:moment_detK},  the corresponding Hamburger moment problem
is determinate and therefore non-vanishing constant functions (and all polynomials) do not belong to the corresponding RKHS. We now investigate the variance of
the discrete BLUE defined in \eqref{eq:estim}, which is given by the ratio of the determinants $H_n$ and~$G_n$.
}

It follows from results in \cite{lau1988extremal} that the determinant of the Hankel matrix defined in \eqref{eq:H_n} has the representation
\begin{equation}\label{L1}
{H}_n=\big|c_{2(i+j)}\big|_{i,j=0}^n=\big|b_{i+j}\big|_{i,j=0}^n= \prod_{i=1}^{n} \left( \tilde{d}_{2i-1} \tilde{d}_{2i} \right)^{n-i+1}\, ,
\end{equation}
where $\tilde{d}_{j}$ are the coefficients of the three-term recurrence relation
\begin{equation}\label{L2}
P_{\ell+1}(t)=(t-\tilde{d}_{2\ell}-\tilde{d}_{2\ell+1})P_{\ell}(t)-\tilde{d}_{2\ell-1}\tilde{d}_{2\ell}P_{\ell-1}(t), \qquad \ell =0,1,\ldots
\end{equation}
of the  monic orthogonal  polynomials
with respect to measure $g(y)dy$ ($\tilde{d}_{0}=0$, $P_0(t)=1$, $P_{-1}(t)=0$).
Observing the three-term recurrence  relation
$$
(\ell+1) L_{\ell+1}^{(\alpha)}(t)=(-t+2\ell+\alpha+1)L_{\ell}^{(\alpha)}(t)-(\ell+\alpha)L_{\ell-1}^{(\alpha)}(t)
$$
for the Laguerre polynomials $L^{(\alpha)}_n(t)$ (orthogonal with respect to $e^{-y}y^{\alpha} dy,$ $ y>0$)
we can identify the coefficients in \eqref{L2}. More precisely, the
 monic polynomials
$$
\overline{L}_{\ell+1}^{(\alpha)}(t)=(-1)^{\ell+1}(\ell+1)!  L_{\ell+1}^{(\alpha)}(t)
$$
 satisfy a three-term recurrence relation of the form \eqref{L2} with $\tilde d_{2k}=k, \tilde d_{2k-1}=k+\alpha$, see \cite{dette1992new}, Lemma 2.2 (b).
As $P_\ell(t)=\overline{L}^{(-1/2)}_\ell (t)$ we have
$
\tilde{d}_{2k}=k,\;\;\tilde{d}_{2k-1}=k-1/2\, ,
$ and therefore obtain
\be \label{H_n}
{H}_{n}\!=\! \prod_{k=1}^{n} \left( k (2k\!-\!1) \right)^{n\!-\!k\!+\!1} \prod_{k=1}^{n} \Big( \frac12 \Big)^{n\!-\!k\!+\!1}=
\Big( \frac12 \Big)^{n(n\!+\!1)/2} \prod_{k=1}^{n} \left( k (2k\!-\!1) \right)^{n\!-\!k\!+\!1}\,.\;\;\;\;\;\;
\ee
Now we move on to the determinant
$
G_{n} = \big|b_{i+j}\big|_{i,j=1}^n  .
$
Note that  we have
$$
b_j=\frac{1}{\sqrt{\pi}} \int_0^\infty y^{j-2} y^{3/2} e^{-y} dy= \frac34 a_{j-2}\,
$$
for $j\geq 2$,
where
$a_k= \int_0^\infty y^k \tilde{g}(y)dy$ and the density $\tilde g_k$ is defined by $ \tilde{g}(y)=\frac{4}{3 \sqrt{\pi}} y^{3/2} e^{-y}$, $y>0$.
Therefore,
$$
G_{n} = \Big( \frac34 \Big)^n \big|a_{i+j}\big|_{i,j=0}^{n-1}= \Big( \frac34 \Big)^n \prod_{l=1}^{n-1} \left(
  \overline{d}_{2l-1} \overline{d}_{2l} \right)^{n-l}\, ,
$$
where
$\overline{d}_{2l-1}=l+3/2$, $ \overline{d}_{2l}=l$.
Consequently,
$$
G_{n} =  \Big( \frac34 \Big)^n  \Big( \frac12 \Big)^{n(n-1)/2} \prod_{k=1}^{n-1} (k(2k+3))^{n-k}
$$
and it follows
\bea
\frac{H_{n}  }{G_{n}}&=&
\Big( \frac43 \Big)^n \Big( \frac12 \Big)^n \,
\Big[  \prod_{k=1}^{n-1}
\frac{(k(2k-1))^{n-k+1}} {
(k(2k+3))^{n-k}} \Big] n(2n-1)\\
&=& \Big( \frac23 \Big)^n n! (2n-1)
 \prod_{k=1}^{n-1}
\frac{(2k-1)^{n-k+1}} {
(2k+3)^{n-k}}\, .
\eea
Since
$$
 \prod_{k=1}^{n-1}
\frac{(2k-1)^{n-k+1}} {
(2k+3)^{n-k}} = \frac{3^n}{(2n-1)(2n+1)!!}
$$
we obtain
\be \label{eq:stein_as}
\frac{H_{n}  }{G_{n}}=  \frac{2^n\, n!}{(2n+1)!!} = \frac{\sqrt{\pi}}{2\sqrt{n}}\Big[1 -\frac{3}{8n} + \frac{25}{128n^2} +O\Big(\frac1{n^3}\Big) \Big], \; n \to \infty.\;\;\;
\ee
The  expansion \eqref{eq:stein_as} details the asymptotic relation formulated as Theorem~3.3 in \cite{xu2017maximum} in the case $p=0$.
Note that   formula  \eqref{eq:stein_as} also  corrects a minor mistake in this reference, which gives $\frac{\sqrt{\pi}}{\sqrt{2n}}$  as the leading term.

\subsection{Spectral measure with Beta distribution }
\label{sec:beta_sp}

\medskip

For measures with a compact support the determinants $H_n$ and $G_n$ can be conveniently evaluated using     the theory of canonical moments, see e.g. \cite{dette1997theory}. Exemplarily, we consider the symmetric Beta $(\alpha,\alpha)$ distribution on the interval $[-1,1]$  with density
\begin{equation}\label{A1}
  \psi^\prime _\alpha (t) = \frac{1}{2^{2 \alpha +1} B(\alpha + 1, \alpha + 1)} (1-t^2)^\alpha,  \qquad -1 < t < 1,
\end{equation}
where  $\alpha > -1$ and
$B(\alpha,\beta)$ denotes the Beta-function. For later purposes we also introduce the
Beta$(\alpha,\beta)$ distribution on the interval $[0,1]$ with density
\be \label{eq:beta1}
\phi_{\alpha,\beta}(t) = \frac1{B(\beta+1,\alpha+1)} t^{\beta}(1-t)^{\alpha}\, , \;\; 0<t<1,
\ee
%here $B(\cdot,\cdot)$ is the Beta-function.
%Moments of this distribution are
%\bea
%m_k(\alpha,\beta)= \int_0^1 t^k \phi_{\alpha,\beta}(t)dt
% = \frac{B(\beta+k+1,\alpha+1)}{B(\beta+1,\alpha+1)} \, .
%\eea
where the  $\alpha, \beta >-1$.
The canonical moments of the Beta-distribution with density \eqref{eq:beta1} are given by
\be
\label{eq:canonical}
p_{2j}=\frac{j}{2j+1+\alpha+\beta}, \;\;
p_{2j-1}=\frac{\beta+j}{2j+\alpha+\beta}\,;
\ee
see e.g. formula (1.3.11) in \cite{dette1997theory}. It is easy to see that the distribution on the interval $[0,1]$ related to the  distribution $\psi_\alpha$  in
\eqref{A1} by the transformation \eqref{h1} is a Beta ($\alpha, - \frac{1}{2}$) distribution. Therefore, it follows from \eqref{eq:canonical} that the corresponding canonical moments are given by

%Beta$(\alpha,\alpha)$ distribution on $[-1,1]$ has density
%$$
%\psi_{\alpha}(t) = \frac1{2^{2\alpha+1}B(\alpha+1,\alpha+1)} (1-t^2)^{\alpha}\, , \;\; -1<t<1.$$
%Moments of this distribution of even order are
%
%\bea
%M_{2k}(\alpha)&=& \int_{-1}^1 t^{2k} \psi_{\alpha}(t)dt
% % \\
%% &=& \frac{1}{2^{2\alpha+1}B(\alpha+1,\alpha+1)}
%% \int_{-1}^1 t^{2k} (1-t^2)^{\alpha}dt\\
%%  &=&
%%  \frac{2}{2^{2\alpha+1}B(\alpha+1,\alpha+1)}
%% \int_{0}^1 t^{2k} (1-t^2)^{\alpha}dt\\
%%   &=&
%% \frac{1}{2^{2\alpha+1}B(\alpha+1,\alpha+1)}   \int_{0}^1 y^{k-1/2} (1-y)^{\alpha}dy
%%\\
%%   &=&
%%\frac{B(1/2,\alpha+1)}{2^{2\alpha+1}B(\alpha+1,\alpha+1)} m_k(-1/2,\alpha)  \\
%%   &=&
%=\frac{B(1/2,\alpha+1)m_k(-1/2,\alpha)}{2^{2\alpha+1}B(\alpha+1,\alpha+1)}   =\frac{B(k+1/2,\alpha+1) }{2^{2\alpha+1}B(\alpha+1,\alpha+1)} \, .
% %=  \frac{\Gamma(k+1/2) \Gamma(2\alpha+2)}{2^{2\alpha+1}\Gamma(\alpha+1) \Gamma(\alpha+k+3/2)}
% %\, ;
%\eea
%%here we have made the substitution $t=\sqrt{y}$, $dt=dy/(2\sqrt{y})$.
%
%
%In view of \eqref{eq:canonical}, canonical moments of Beta$(-1/2,\alpha)$ distribution on $[0,1]$ are
\begin{equation}
\label{A2}
p_{2j}=\frac{j}{2j+1/2+\alpha}, \;\;
p_{2j-1}=\frac{j-1/2}{2j-1/2 +\alpha}~.
\end{equation}
Now Theorem 1.4.10 in \cite{dette1997theory} gives
\begin{equation} \label{A3}
H_n  =   |(b_{i+j})^n_{i,j=0}| =
 \prod_{i=1}^{n} \left( q_{2i-2}p_{2i-1}q_{2i-1}p_{2i}\right)^{n+1-i},
 \end{equation}
 where $q_0=1, q_j=1-p_j$ $ (j \geq 1)$ and (observing \eqref{A2})
 \begin{equation}\label{A6a}
q_{2i-2}p_{2i-1}q_{2i-1}p_{2i}= {\frac {4i \left( i+\alpha \right)  \left( 2\,i-1+2\,\alpha \right)
 \left( 2\,i-1 \right) }{ \left( 4\,i+1+2\,\alpha \right)  \left( 4\,i
-1+2\,\alpha \right) ^{2} \left( 4\,i-3+2\,\alpha \right) }}\, ,\;\; i=1,2 \ldots~
\end{equation}

 For the calculation of the determinant $G_n=|(b_{i+j})^n_{i,j=1}|$ we note  the relation
 \begin{equation}\label{A4}
   b_i = \frac{B(\frac{5}{2}, \alpha +1)}{B(\frac{1}{2}, \alpha + 1)} \tilde b_{i-2} \qquad i=2,3,\ldots
 \end{equation}
 where $\tilde b_0, \tilde b_1, \ldots$ are the moments of the {Beta($\alpha, {3}/{2}$)} distribution. Consequently, we obtain from Theorem 1.4.10 in \cite{dette1997theory} that
  \begin{eqnarray} \label{A5}\nonumber
 G_n &=& |(b_{i+j})^n_{i,j=1}| = |(\tilde b_{i+j})^{n-1}_{i,j=0}| = \left[\frac{B(\frac{5}{2},\alpha+1)}{B(\frac{1}{2},\alpha +1) }\right]^{n}
  \times \prod_{i=1}^{n-1} \left( \tilde q_{2i-2} \tilde p_{2i-1} \tilde q_{2i-1} \tilde p_{2i}\right)^{n-i}\\
  &=& \left[ \frac{3}{(2 \alpha +3)(2 \alpha +5)} \right]^{n} \times \prod_{i=2}^{n} \left( \tilde q_{2i-4} \tilde p_{2i-3} \tilde q_{2i-3} \tilde p_{2i-2}\right)^{n+1-i}\, ,
\end{eqnarray}

%with canonical moments of Beta$(1/2,\alpha)$ distribution
%\bea
%\label{eq:canonical5}
%\tilde p_{2j}=\frac{j}{2j+3/2+\alpha}, \;\;
%\tilde p_{2j-1}=\frac{\alpha+j}{2j+1/2 +\alpha}
%\eea
%and
%$$
%\tilde q_{2i-2} \tilde p_{2i-1}\tilde q_{2i-1}\tilde p_{2i}= {\frac {4i \left( i+\alpha \right)  \left( 2\,i+1+2\,\alpha \right)
% \left( 2\,i+1 \right) }{ \left( 4\,i+3+2\,\alpha \right)  \left( 4\,i
%+1+2\,\alpha \right) ^{2} \left( 4\,i-1+2\,\alpha \right) }}
%\, ,\;\; i=1,2 \ldots
%$$
%
%

where $\tilde p_1, \tilde p_2$ are the canonical moments of Beta$(\alpha ,3/2)$ distribution, that is
\bea
\label{eq:canonical5}
\tilde p_{2i}=\frac{j}{2i+5/2+\alpha}\,, \;\;
\tilde p_{2i-1}=\frac{3/2 +i}{2i+3/2 +\alpha}\, ,
\eea
and
\begin{equation}\label{A6}
\tilde q_{2i-2} \tilde p_{2i-1}\tilde q_{2i-1}\tilde p_{2i}=
{\frac {4i \left( i+\alpha \right)  \left( 2\,i+3+2\,\alpha \right)
 \left( 2\,i+3 \right) }{ \left( 4\,i+5+2\,\alpha \right)  \left( 4\,i
+3+2\,\alpha \right) ^{2} \left( 4\,i+1+2\,\alpha \right) }}
\, ,\;\; i=1,2 \ldots
\end{equation}

%Let
%\bea
%c_{n, \alpha}&=&
%\frac{B(1/2,\alpha+1)}{2^{2\alpha+1}B(\alpha+1,\alpha+1)}
% \left[\frac{B(1/2,\alpha+1)}{ B(5/2,\alpha+1)}\right]^{n} \\
% % &=& \frac{\Gamma(1/2)  \Gamma(2\alpha+2)}{2^{2\alpha+1}\Gamma(3/2+\alpha) \Gamma(\alpha+1)}
%%  \left[\frac{\Gamma(7/2+\alpha) \Gamma(1/2) \Gamma(\alpha+1)}{ \Gamma(3/2+\alpha) \Gamma(5/2) \Gamma(\alpha+1) }\right]^{n}\\
% &=& \frac{\sqrt{\pi}  \Gamma(2\alpha+2)}{2^{2\alpha+1}\Gamma(3/2+\alpha) \Gamma(\alpha+1)}
%  \left[\frac{(3+2\alpha)(5+2\alpha) }{ 3  }\right]^{n}\, .
%\eea
Consequently, it follows from \eqref{A6},   \eqref{A5} and \eqref{A6a}
\bea
\!\!\frac{H_n}{G_n}&=&  \Big [\frac{(2 \alpha +3)(2 \alpha +5)}{3} \Big ]^n \left( q_{0}p_{1}q_{1}p_{2}\right)^n \prod_{i=2}^{n} \Big[\frac{ q_{2i-2}p_{2i-1}q_{2i-1}p_{2i} } { \tilde q_{2i-4} \tilde p_{2i-3} \tilde q_{2i-3} \tilde p_{2i-2} }\Big]^{n+1-i}\\
\!\!\!\!\!\!&=& \!\!
\Big[{\frac {4(1\!+\!\alpha)}{ 3  \left( 3\!+\!2\,\alpha
 \right) ^{2}}} \Big]^n
 \prod_{i=2}^{n}
 \Big[
{\frac {i \left( i\!+\!\alpha \right) \left( i\!-\!1/2 \right) \left( i\!+\!\alpha\!-\!1/2 \right)
  }{ \left( i\!-\!1 \right)  \left( i\!-\!1\!+\!\alpha
 \right)  \left( i\!+\!1/2 \right) \left( i\!+\!\alpha\!+\!1/2 \right)  }}
  \Big]^{n+1-i}\, .
\eea \nopagebreak
Observing the relations

%\bea
%\frac{H_n}{G_n}&=& \left( q_{0}p_{1}q_{1}p_{2}\right)^n \prod_{i=2}^{n} \left[\frac{ q_{2i-2}p_{2i-1}q_{2i-1}p_{2i} } { \tilde q_{2i-4} \tilde p_{2i-3} \tilde q_{2i-3} \tilde p_{2i-2} }\right]^{n+1-i}\\
%&=&
%\left[{\frac {4(1+\alpha)}{ \left( 5+2\,\alpha \right)  \left( 3+2\,\alpha
% \right) ^{2}}} \right]^{n}
% \prod_{i=2}^{n}
% \left[       1+ {\frac { \left( 4\alpha^2-1 \right)
% \left( \alpha+2\,i-1 \right) }{
% \left( (4i+2\,\alpha)^2-1 \right)  \left( i-1 \right)  \left( i-1+
%\alpha \right) }}
%  \right]^{n+1-i}
%\eea
\bea
 \prod_{i=2}^{n}
 \Big[
{\frac {i  }{  i-1
}}
  \Big]^{n+1-i}&=&n!\,,\\
   \prod_{i=2}^{n}
 \Big[
{\frac {i-1/2  }{  i+1/2
}}
  \Big]^{n+1-i}&=&\frac{3^n}{(2n+1)!!}\, ,\\
   \prod_{i=2}^{n}
 \Big[
{\frac { i+\alpha  }{ i-1+\alpha
}}
  \Big]^{n+1-i}
  &=&\frac{\Gamma(n+1+\alpha)}{(1+\alpha)^n\Gamma(1+\alpha)}\, ,\\
   \prod_{i=2}^{n}
 \Big[
{\frac { i+\alpha -1/2 }{ i+\alpha+1/2
}}
  \Big]^{n+1-i}
  &=&\frac{(3+2\alpha)^{n}\Gamma  (3/2+\alpha)}{2^n \Gamma  (n+3/2+\alpha)}\, .
  \eea
we obtain
\be \nonumber
\frac{H_n}{G_n}&=&
\Big[{\frac {4(1\!+\!\alpha)}{ 3  \left( 3\!+\!2\,\alpha
 \right)  }} \Big]^{n}
  \frac{n! 3^n \Gamma(n\!+\!1\!+\!\alpha) (3\!+\!2\alpha)^{n}\Gamma  (3/2\!+\!\alpha)}{(2n\!+\!1)!! (1\!+\!\alpha)^n\Gamma(1\!+\!\alpha)2^n \Gamma  (n\!+\!3/2\!+\!\alpha)}
 % \frac{\Gamma(n+1+\alpha)}{(1+\alpha)^n\Gamma(1+\alpha)}
%\frac{(3+2\alpha)^{n}\Gamma  (3/2+\alpha)}{2^n \Gamma  (n+3/2+\alpha)}
%{\rm const}
%  \left[\frac{(3+2\alpha)(5+2\alpha) }{ 3  }\right]^{n}\left[{\frac {4(1+\alpha)}{ \left( 5+2\,\alpha \right)  \left( 3+2\,\alpha
% \right) ^{2}}} \right]^{n}  \frac{n!3^n(3/2+\alpha)^{n} \Gamma(n+1+\alpha)}{(2n+1)!!(1+\alpha)^n \Gamma  (n+3/2+\alpha)}
 \\ \nonumber
     &=& \frac{\sqrt{\pi} } {2^{2\alpha+1} B (\alpha+1,\alpha+1)}
   \times   \frac{(2n)!! }{(2n+1)!!} \cdot \frac{ \Gamma(n+1+\alpha)}{\Gamma  (n+3/2+\alpha)} \\
     &=&\frac{{\pi} } {2^{2\alpha+2} B (\alpha+1,\alpha+1)} \times \frac{1}{n} \Big(1 +O\Big(\frac1n \Big)\Big),\; n \to \infty \,   \label{eq:sa_beta}
\ee
where the expansion in the last line follows by straightforward but tedious calculation using Stirling's formula.

We finally mention the special cases  $\alpha=0$ (the spectral measure is a uniform spectral density on the interval $[-1,1]$ with corresponding kernel function $k(x)=\sin(x)/x   $)
and
 $\alpha=-1/2$   (the spectral measure is the  arcsine  distribution on $[-1,1]$ and the corresponding kernel   is $k(x)=2J_1(x)/x$, where $J_{\alpha}(\cdot)$ is the Bessel function of the first kind)
for which the  expansions are  given, respectively,  by
\be \label{eq:stein_asU}
\frac{H_n}{G_n} &=& \left[ \frac{(2 n)!!}{(2n+1)!!}\right]^2=
 \frac{{\pi}}{4{n}}\! +\! O\Big(\frac1{n^2}\Big) ,
\\
\nonumber
\frac{H_n}{G_n} &=&
 \Big(\frac83\Big)^n \Big[\frac {1}{ 8} \Big]^{n} n! \frac{ \,3^{n}}{(2n+1)!! }\frac{2^n \Gamma(n+1/2) }{\Gamma(1/2)} \frac{1}{n!}
  %\\  &=&
  =  \frac{1}{2n+1}
  %\label{eq:stein_asU12}
 = {1 \over 2n} +\! O\Big(\frac1{n^2}\Big)
\ee
as $ n \to \infty. $
Interestingly, the ratio ${H_n}/{G_n}$ in \eqref{eq:stein_asU} is the squared  ratio ${H_n}/{G_n}$ of \eqref{eq:stein_as}.

\section{Some results of numerical studies and discussions}
\label{sec:numer}
%
%
%{Note that under the assumptions of Theorems~\ref{th:one_dim} and~\ref{th:one_dim1} all functions of the form  $f\!+\!g $ with  $g \in H(K)$ and
%$f $ being either a non-zero constant or a polynomial do not belong to $ H(K)$. }
%
%
%
%
%\bigskip

 {We have made extensive numerical studies to assess the uncertainty quantification in GP regression, as introduced in Section~\ref{sec:int_2}
 for functions $f\in H(K)$ and $f\notin H(K)$; some of our results are illustrated in the figures below. At the end of this section, we summarize our conclusions.
Different  kernels $K$ (including Mat\'{e}rn kernels and the kernels discussed at the end of Section~\ref{sec:suff}) have been investigated as well.
 In the
figures  below, we use $X=[0,1]$, Gaussian and Cauchy kernels (see (E.1) and (E.2) at the end of Section~\ref{sec:suff}) and the  two functions
 $$
 f_1(x)=\exp\{- 2 (x- 1/3 )^2\} ~~\text{ and }  ~~f_2(x)=1- 2 (x- 1/3 )^2.
 $$
These two functions look similar but  we note that $f_1\in H(K)$ and $f_2\notin H(K)$ for both  kernels with the correlation lengths  considered below.
% (to prove that $f_1\in H(K)$ in the case of  the
%Gaussian kernel  we can use Theorem 3 in \cite{minh2010some}). 
Visually, the chosen kernels also look similar but it turns out that they exhibit completely different behaviour.}

\begin{figure}[h]
\centering

  \includegraphics[width=0.49\linewidth]{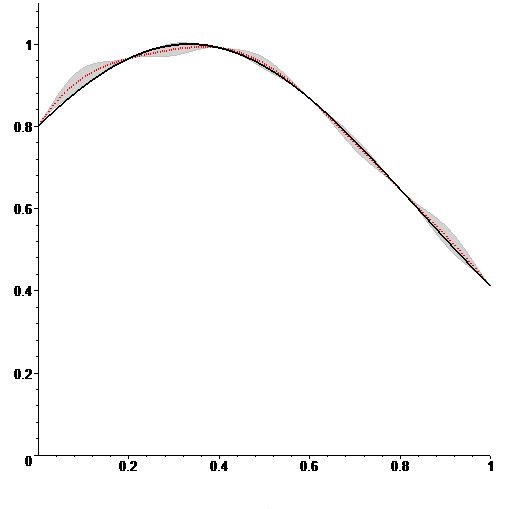}
 \includegraphics[width=0.49\linewidth]{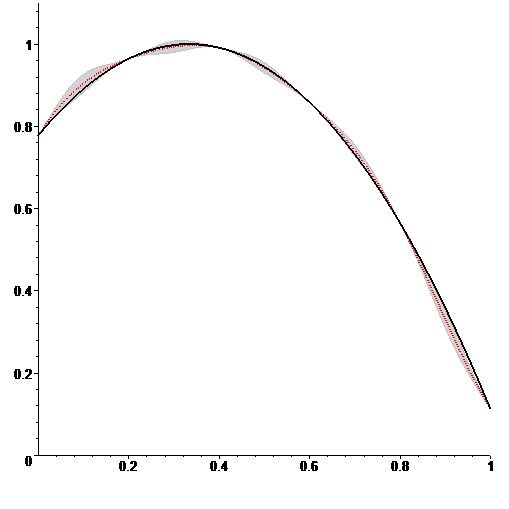}
 \caption{\it Kriging confidence regions for kernel  approximation of $f_1$ (left) and $f_2$ (right): Gaussian kernel, $\lambda=15$, $N=6$.}
\label{error_f_Gaus_6}

\end{figure}

\begin{figure}[h]
\centering

  \includegraphics[width=0.49\linewidth]{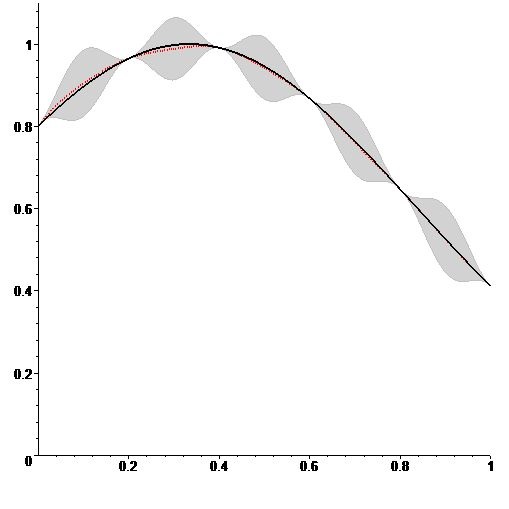}
 \includegraphics[width=0.49\linewidth]{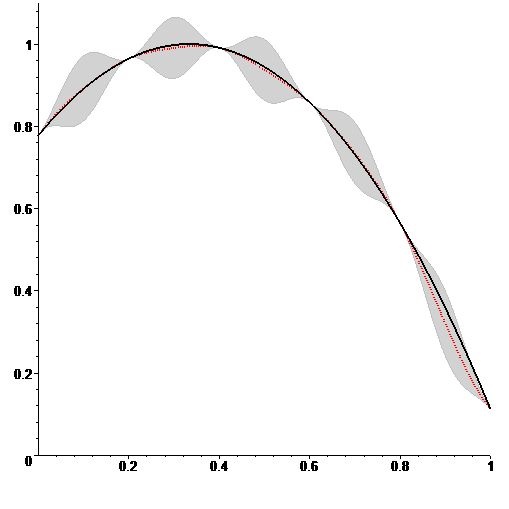}
 \caption{\it Kriging confidence regions for kernel  approximation of $f_1$ (left) and $f_2$ (right): Cauchy kernel, $\lambda=20$, $N=6$.}
\label{error_f_Cauchy_6}

\end{figure}

\begin{figure}[h]
\centering

  \includegraphics[width=0.49\linewidth]{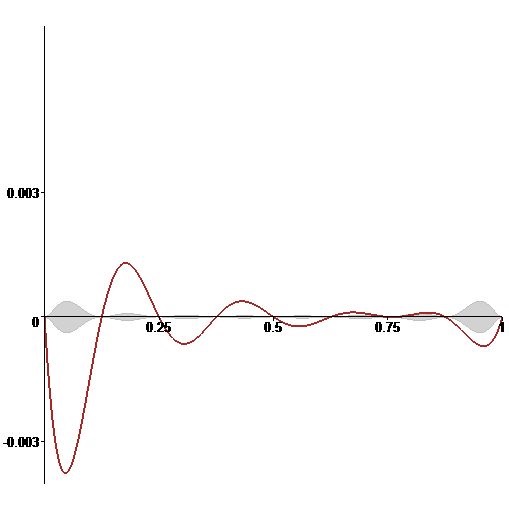}
  \includegraphics[width=0.49\linewidth]{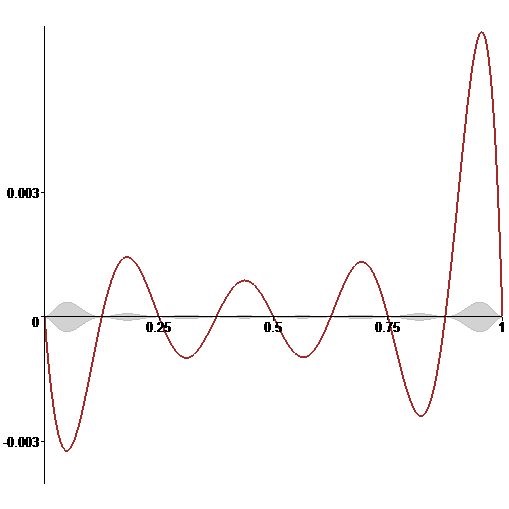}
 \caption{\it Deviation  and confidence   bounds  in kernel approximation:   Gaussian kernel, $\lambda=15$,  $N=9;$ $f_1$ (left) and $f_2$ (right)}
\label{error_Gaus_9}

\end{figure}

\begin{figure}[h]
\centering

  \includegraphics[width=0.49\linewidth]{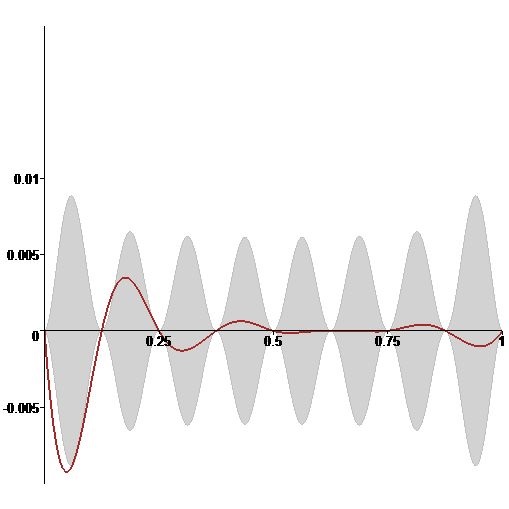}
  \includegraphics[width=0.49\linewidth]{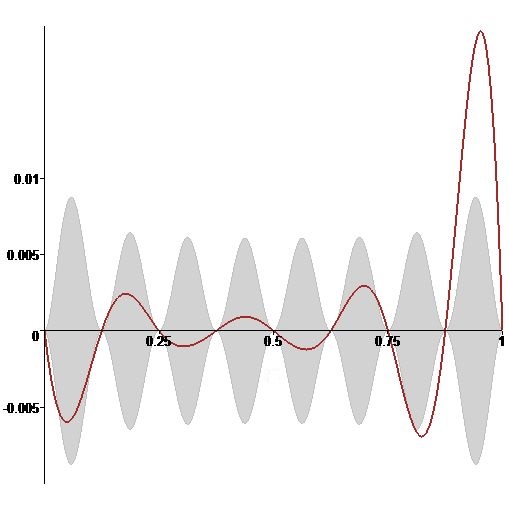}
 \caption{\it Deviation  and confidence   bounds  in kernel approximation:   Cauchy kernel, $\lambda=20$,  $N=9;$ $f_1$ (left) and $f_2$ (right)}
\label{error_Cauchy_9}

\end{figure}

\begin{figure}[h]
\centering
  \includegraphics[width=0.49\linewidth]{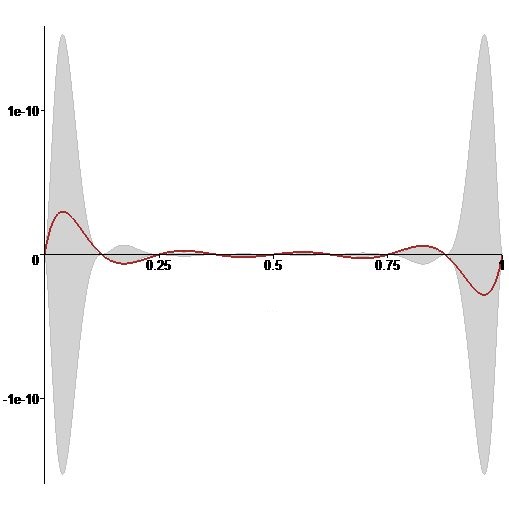}
  \includegraphics[width=0.49\linewidth]{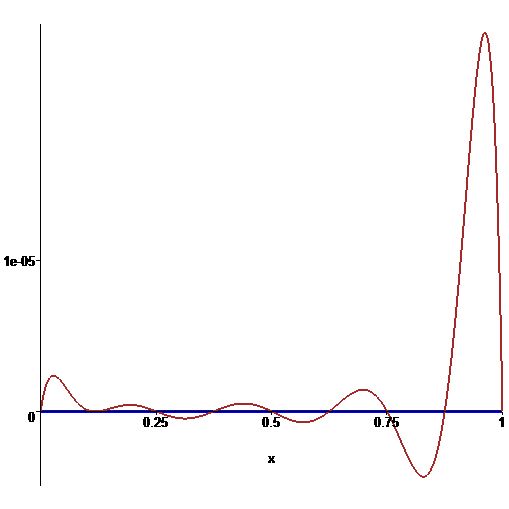}
 \caption{\it Deviation  and confidence   bounds  in kernel approximation:   Gaussian kernel, $\lambda=2$,  $N=9;$ $f_1$ (left) and $f_2$ (right).}
\label{error_Gaus_9_2}

\end{figure}

 {In Figs.~\ref{error_f_Gaus_6} and \ref{error_f_Cauchy_6} we plot either $f_1$  or $f_2$  in solid black, the kernel approximation $\mu_N(x)$ computed by \eqref{det1} in dotted red
 and the so-called kriging confidence regions
$\mu_N(x) \pm 3\widehat{\sigma^2_N}\, C_N(x,x)$ in grey, where  $C_N(x,x)$ is the kernel variance computed by  \eqref{det11} and $\widehat{\sigma^2_N}$ is the MLE of $\sigma^2$ computed by  \eqref{1s}. The
main reason for providing Figs.~\ref{error_f_Gaus_6} and \ref{error_f_Cauchy_6} is the demonstration of the big difference in the width of the confidence regions for the Gaussian and Cauchy kernels.
In Figs.~\ref{error_Gaus_9}, \ref{error_Cauchy_9} and \ref{error_Gaus_9_2}, we plot the deviation $f(x)-\mu_N(x)$ in brown and confidence bounds $\mu_N(x)-f(x) \pm 3\widehat{\sigma^2_N}\, C_N(x,x)$ in filled grey.
Again, the left and right panels show the results for the functions $f_1$ and $f_2$, respectively.
The points $x_j$, where observations of $f$ are taken, are equally spaced on the interval  $[0,1]$ with $x_j=(j-1)/(N-1)$, $j=1, \ldots, N$.
}

 {
The results for the Gaussian kernel $K(x,y)=\sigma^2\exp \{-15(x-y)^2\} $ are depicted on Figs.~\ref{error_f_Gaus_6} and~\ref{error_Gaus_9}.
The corresponding  results  for the   Cauchy kernel $K(x,y)=\sigma^2/(1+20 (x-y)^2) $ can be found
Figs.~\ref{error_f_Cauchy_6} and \ref{error_Cauchy_9}. It is clear from comparing left and right panels in Figs.~\ref{error_f_Gaus_6}--\ref{error_Cauchy_9} that kernel approximations for $f_1\in H(K)$ are significantly more accurate than for $f_2\notin H(K)$. The two chosen kernels (Gaussian with $\lambda=15$ and Cauchy with $\lambda=20$) look very similar but have different tail behaviour of the corresponding spectral density: the tail of the spectral density of the Cauchy kernel  has a heavier tail.
The confidence regions for the regression function  constructed by the Cauchy kernel are rather wide and resemble the regions for the Mat\'{e}rn kernels with shape parameters 3/2 and 5/2 having similar correlation lengths. The  confidence regions in the case of  Gaussian kernel are much  narrower (in fact, far too narrow) and the   confidence regions for the kernels in (E.3) and (E.4) of Section~\ref{sec:suff} are even narrower; the spectral measures for these kernels have finite support.
}

 {In Figure~\ref{error_Gaus_9_2} we plot the deviations and confidence regions for kernel approximation with Gaussian kernel $K(x,y)=\sigma^2\exp \{-2(x-y)^2\} $; the Gaussian kernel with $\lambda=2$ is perfectly suited for function
$f_1$. A naive
visual inspection of the two functions might suggest the Gaussian kernel
should also work well for $f_2$ but,  as we can observe from Figure~\ref{error_Gaus_9_2} (right), it is not so for
$f_2 \notin H(K)$. The confidence region on Figure~\ref{error_Gaus_9_2} (right) cannot be seen as  the deviations $|\mu_N(x)-f_2(x)|$ are on average $10^5$ times larger than $3\widehat{\sigma^2_N}\, C_N(x,x)$.
}

 {From the numerical studies partially illustrated in these figures we make the following conclusions concerning uncertainty quantification in GP regression models with infinitely differentiable translation-invariant kernels:
\begin{itemize}
\item if $f \notin H(K)$, then the kriging confidence regions for $f$ are always  inaccurate;
  \item the heavier are the tails of the spectral measure of the kernel, the wider are the confidence regions;
  \item if the  tail of the spectral measure of the kernel is  light  and the function $f$ does not belong to the respective RKHS, then the kernel approximation of $f$ appears to be rather inaccurate and the confidence regions  seem to be missing $f$  almost entirely;
  \item if the function $f$ does not belong to the respective RKHS, then $\widehat{\sigma^2_N}\to \infty$ as $N\to \infty$,
   but this has little effect on the size of the confidence intervals, at least  for small $N$;
  \item for kernels with light tails of the respective spectral measures, the kernel approximation is accurate and confidence regions are adequate only if the shape   of $f$ precisely matches the shape of the kernel functions $K(x,\cdot)$, as in Figure~\ref{error_Gaus_9_2} (left).
%\item
%smoothness  of the kernel and  the inclusion $f \in H(K)$ have little effect on  uncertainty quantification in GP regression in the  case when there is measurement error (i.e., a white noise process representing nugget effect).
\end{itemize}
}

{\bf Acknowledgements}
This research was partially supported by the research grant DFG DE 502/27-1 of the German Research Foundation (DFG).
{The authors  are  grateful to Timo Karvonen (Alan Turing Institute) for intelligent discussions, spotting an essential typo and pointing out several important references.
We would also like to thank both referees  and  {especially} the associate editor for their constructive and very valuable comments on  earlier versions of this paper.
}

\end{document}